\documentclass[letterpaper,onecolumn,12pt]{article}    
\usepackage{amssymb}
\usepackage{amsmath}
\usepackage{graphicx}  
\usepackage{epstopdf}
\usepackage{hyperref}
\usepackage{color}
\usepackage{bbold}
\usepackage{units}
\usepackage{a4wide}

\newtheorem{theorem}{Theorem}
\newtheorem{corollary}{Corollary}
\newtheorem{lemma}{Lemma}
\newtheorem{remark}{Remark}
\newtheorem{assumption}{Assumption}
\newtheorem{proposition}{Proposition}
\newtheorem{definition}{Definition}
\newtheorem{example}{Example}

\definecolor{changecolor}{rgb}{0.0, 0.0, 0.0}

\newcommand{\defeq}{\overset{\mathrm{def}}{=}}
\newcommand{\tr}{\intercal}

\newcommand{\rank}{\mathrm{rank}}
\newcommand{\qed}{$\diamond$}

\newcommand{\change}[1]{{\color{changecolor}#1}}

\title{\textbf{Configuration-Constrained Tube MPC}}

\author{Mario E.~Villanueva$^{a,c}$, Matthias A.~M\"uller$^b$, Boris Houska$^c$}
\date{\small $^a$ IMT Lucca \\[0.1cm] $^b$ Leibniz University Hannover \\[0.1cm] $^c$ ShanghaiTech University}

\begin{document}

\addtolength{\textheight}{-2cm}

\maketitle


\begin{abstract}
This paper is about robust Model Predictive Control (MPC) for linear systems \change{with additive and multiplicative uncertainty}. \change{A novel class of configuration-constrained polytopic robust forward invariant tubes is introduced, which admit a joint parameterization of their facets and vertices. They are the foundation for the development of novel} \change{Configuration-Constrained T}ube MPC \change{(CCTMPC)} controllers \change{that freely optimize the shape of their polytopic tube, subject to conic vertex configuration constraints, as well as associated vertex control laws by solving convex optimization problems online.} \change{It is shown that CCTMPC is---under appropriate assumptions---systematically less conservative than Rigid- and Homothetic- Tube MPC. Additionally, it is proven that there exist control systems for which CCTMPC is less conservative than Elastic Tube MPC, Disturbance Affine Feedback MPC, and Fully Parameterized Tube MPC.}\\[-0.3cm]
\end{abstract}

\section{Introduction}
\label{sec:introduction}
During the last two decades, Tube MPC~\cite{Langson2004} has emerged as a sensible alternative to robust dynamic programming~\cite{Bertsekas2012,Diehl2004}
and min-max feedback MPC~\cite{Goulart2006} for formulating, analyzing and approximating the robust control synthesis problem~\cite{Scokaert1998}. It is a set-based framework~\cite{Blanchini2008}, whose underlying principle consists of replacing system trajectories by so-called robust forward invariant tubes. These are set-valued functions in the state-space, enclosing all future states of the system, for a given feedback law, independently of the uncertainty realization~\cite{Mayne2005}.

\bigskip
\noindent
Practical Tube MPC formulations rely on the parameterization of their set-valued tubes and, in the context of many existing approaches, also the feedback law. Thus, Tube MPC is harder to formulate and solve than certainty-equivalent MPC. As such, from a practical perspective, one has to ask whether the effort of investing into a robust MPC formulation pays out in the first place~\cite{Mayne2015}. However, at least for linear systems, a variety of tractable convex
Tube MPC formulations exists. This includes the early tube MPC formulation from~\cite{Chisci2001} as well as the so-called \change{Rigid-Tube MPC (RTMPC)~\cite{Mayne2005}, Homothetic Tube MPC (HTMPC)~\cite{Rakovic2012,Rakovic2013}, Elastic Tube MPC (ETMPC)~\cite{Rakovic2016}, and Fully Parameterized Tube MPC (FPTMPC)~\cite{Rakovic2011FP,Rakovic2012a}} formulations.

\bigskip
\noindent
In general, it is impossible to optimize over arbitrary sets and feedback policies. Consequently, there is no unique way of answering what is the best way of implementing robust MPC. Instead, one has to trade-off between introducing conservatism 
and improving computational run-time performance~\cite{Kohler2019}. For instance, \change{ETMPC} is less conservative than \change{RTMPC}, but comes at the cost of introducing more optimization variables. Similarly, affine feedback policies \change{are, in general,} sub-optimal~\cite{Bertsimas2012}, but many robust MPC schemes use affine ancillary feedback \change{\mbox{laws---or} affine disturbance feedback laws as in the context of Disturbance Affine Feedback MPC \mbox{(DAFMPC)~\cite{Goulart2006}---in}} order to arrive at a tractable reformulation. And, last but not least, the parameterization of the tube is often based on polytopes~\cite{Rakovic2012}, ellipsoids~\cite{Villanueva2017}, or other classes of computer representable sets, but the choice of such set parameterizations affects the run-time and performance of the associated controller.

\paragraph*{Main Contribution.}

\change{The main contribution of this paper is a novel class of Configuration-Constrained Tube MPC (CCTMPC) controllers for linear discrete-time systems with additive and multiplicative uncertainty that admit an exact reformulation as convex optimization problem while avoiding a direct parameterization of the feedback law. The computational complexity of CCTMPC scales linearly with respect to the length of the controller's prediction horizon while its level of conservatism solely depends on the number of facet directions and the vertex configuration of the polytopes that are used to parameterize the tube. As we will prove in this paper, CCTMPC \mbox{is---under} suitable assumptions that, however, merely aim at making these controllers comparable at \mbox{all---never} more conservative than RTMPC and HTMPC. Here, we say that a given Tube MPC scheme ``$\mathsf{A}$'' is never more conservative than another Tube MPC scheme ``$\mathsf{B}$'' only if we are capable to prove that the set of all possible tubes and all possible feedback laws that are representable by scheme $\mathsf{A}$ contains the corresponding set of representable tubes and feedback laws of scheme $\mathsf{B}$. And it is also in this rather strict sense that we shall prove that CCTMPC is never more conservative than ETMPC for systems with two states. It needs to be also stated, however, that in this particularly strict sense, one cannot compare CCTMPC and ETMPC in higher dimensional state spaces. Nevertheless, we will discuss numerical examples, which indicate that CCTMPC performs better than ETMPC for a couple of selected problems in higher dimensions and with respect to a selected objective. Besides, a principal advantage of CCTMPC is that it can directly tackle both additive as well as multiplicative uncertainties. In order to facilitate similar extensions of HTMPC and ETMPC for systems with multiplicative uncertainty one would first need to find an affine control law that is robust with respect to the mentioned additive and multiplicative uncertainties, which can be a difficult task. This is in contrast to CCTMPC, which is not based on the availability of such robust affine control laws. Moreover, as we shall establish in this paper, too, there exist linear control systems for which CCTMPC can be proven to be both strictly less conservative and strictly less computationally demanding than DAFMPC and FPTMPC.}

\paragraph*{Overview.}
This paper \change{is organized as follows.}

\begin{itemize}


\item Section~\ref{sec::tmpc} introduces the problem formulation.


\item Section~\ref{sec::polyhedra} reviews the definition of template polyhedra.
In this context, our contribution is the introduction of novel conic configuration domains, which correspond to the set of parameters of template polyhedra that share a given partial face configuration. Theorem~\ref{thm::regularTemplates} provides a unique characterization of such configuration domains. Moreover, \change{Theorem~\ref{thm::convhull} elaborates on the use of such configuration domains for computing vertex representations of polytopes that are given in half-space representation.}

\item Section~\ref{sec::Tube} \change{presents an exact and computationally tractable condition under which a sequence of configuration-constrained polytopes form a robust forward invariant tube. This relies on the use of a vertex control law; see Corollary~\ref{cor::F}.}

\item Section~\ref{sec::MPCsummary} \change{discusses the practical implementation of CCTMPC by using convex optimization.} The fact that this controller is recursively feasible and asymptotically stable is established in Theorem~\ref{thm::stability}.

\item Section~\ref{sec::caseStudy} \change{presents an in-depth discussion of three tutorial examples, which are additionally supported by numerical illustrations. In detail, it is explained and visualized under which assumptions CCTMPC is less conservative than RTMPC, HTMPC, and ETMPC. Moreover, an explicit example for an uncertain linear system with $4$ states and $1$ control input is constructed, for which it can be shown that CCTMPC is less conservative than both DAFMPC and FPTMPC.}

\end{itemize}

\noindent
Finally, Section~\ref{sec::conclusions} concludes the paper.

\section{Tube MPC\label{sec::tmpc}}
This section introduces the basic notation for uncertain linear systems, reviews the concept of robust forward invariance, and discusses our problem formulation.

\subsection{Uncertain Linear Systems}
\label{sec::LinearSystem}
This paper is concerned with systems of the form
\begin{equation} 
\label{eq::system}
x_{k+1} = A\change{_k} x_k + B\change{_k} u_k + C w_k \;.
\end{equation}
Here, $x_k \in \mathbb X$ denotes the state at time $k \in \mathbb N$ and $u_k \in \mathbb U$ the control input. The state and control constraint sets, $\mathbb X \subseteq \mathbb R^{n_x}$ \change{and} $\mathbb U \subseteq \mathbb R^{n_u}$, \change{are assumed to be closed and convex.} \change{In the most general setting, the matrices
$$[A_k,B_k] \in \Delta \subseteq \mathbb R^{n_x \times (n_x+n_u)}$$
and vectors $w_k \in \mathbb W$ are both unknown. Here, $\Delta $ denotes a matrix polytope,
\[
\Delta \ \defeq \ \mathrm{conv}\left( \ [\overline A_1, \overline B_1], \ [\overline A_2, \overline B_2], \ \ldots, \ [\overline A_l, \overline B_l] \ \right),
\]
with given vertices, $[\overline A_1, \overline B_1], \ [\overline A_2, \overline B_2], \ \ldots, [\overline A_l, \overline B_l]$.} The set $\mathbb W \subseteq \mathbb R^{n_w}$ \change{is assumed to be compact. Additionally, $C\in\mathbb{R}^{n_x\times n_w}$ is assumed to be given and constant}.

\subsection{Robust Forward Invariant Tubes}
\label{sec::rfit}
\change{Let $\mathcal U = \{ \mu \mid \mu: \mathbb R^{n_x} \to \mathbb U \}$ denote the set of control laws; that is, the set of maps from $\mathbb R^{n_x}$ to $\mathbb U$. In the following, we use the notation
\[
f(X,\mu) \ \defeq \ \left\{
\ A x + B\mu(x) + C w \ \middle|
\begin{array}{l}
\ x \in X, \ w \in \mathbb W \\[0.16cm]
\left[ A, B \right] \in \Delta 
\end{array}
\right\}
\]
to denote the closed-loop set propagation function, which is defined for all sets $X \subseteq \mathbb R^{n_x}$ and all control laws $\mu \in \mathcal U$. The following definitions are standard in the set-theoretic control literature~\cite{Blanchini1999}.}

\begin{definition}\label{def::rci}
\textit{A set $X_\mathrm{s} \subseteq \mathbb R^{n_x}$ is called a robust control invariant (RCI) set of~\eqref{eq::system} if \change{there exists a $\mu_\mathrm{s} \in \mathcal U$ for which $X_\mathrm{s} \supseteq f(X_\mathrm{s},\mu_\mathrm{s})$}.}
\end{definition}

\begin{definition}\label{def::rfit}
\textit{A sequence of sets $X_0, X_1, \ldots X_N \subseteq \mathbb R^{n_x}$, with $N \in \mathbb N \cup \{ \infty \}$, is called a robust forward invariant tube of~\eqref{eq::system} if \change{there exists a sequence of control laws $\mu_0, \mu_1, \ldots \mu_N \subseteq \mathcal U$ such that $X_{k+1} \supseteq f(X_k,\mu_k)$ for all indices $k \in\{ 0,1, \ldots, N-1 \}$.}}
\end{definition}

\noindent
The function $f$ should not be mixed up with the function
\[
F(X) \defeq \left\{
\; X^+ \subseteq \mathbb R^{n_x} \; \middle|
\begin{array}{l}
\forall x \in X, \; \exists u \in \mathbb U: \\[0.16cm]
\change{\forall [A,B] \in \Delta,} \; \forall w \in \mathbb W, \\[0.16cm]
A x + Bu + C w \in X^+
\end{array}
\right\} \, ,
\]
which \change{is also defined for all} sets $X \subseteq \mathbb R^{n_x}$. \change{The functions $f$ and $F$ are, however, closely related, since for two given sets $X,X^+ \subseteq \mathbb R^{n_x}$ we have
\[
\exists \mu \in \mathcal U: \ X^+ \supseteq f(X,\mu) \quad \Longleftrightarrow \quad X^+ \in F(X) \; .
\]
Notice that while $f$ is set-valued, the values of $F$ are sets-of-sets. Despite this 
apparent complication, working with $F$ is sometimes more elegant than working with $f$.
For instance, instead of Definition~\ref{def::rci}, one could also say that $X_\mathrm{s} \subseteq \mathbb R^{n_x}$ is an RCI set if \mbox{$X_\mathrm{s} \in F(X_\mathrm{s})$}. This is not only shorter than the sentence from Definition~\ref{def::rci}, but it also avoids to explicitly introduce a control law. In fact, as we shall see later on in this paper, the switch from $f$ to $F$ is more than a switch of notation. Namely, it is intended to highlight the fact that, in the context of Tube MPC, a direct parameterization of control laws can eventually be avoided as long as one is able to find a computationally tractable representation of $F$. 
}

\subsection{Tube MPC}
\label{sec::TubeMPC}
The focus of this paper is on formulating, analyzing, and solving Tube MPC problems of the form 
\begin{align}
\label{eq::tocp}
\hspace{-0.2cm}
\begin{array}{cl}
\underset{\change{X,\mu}}{\min} & \change{L_0(X_0) + \displaystyle\sum\limits^{N-1}_{k=0} L(X_k,\mu_k) + L_N(X_N)}
\\[0.7cm]
\text{s.t.} &
\left\{ 
\begin{array}{l}
\forall k\in\{ 0, 1, \ldots, N-1\}, \\[0.16cm]
\change{X_{k+1} \supseteq f(X_k,\mu_k), \ \hat x \in X_0, \ X_k\subseteq \mathbb X} , \\[0.16cm]
\change{X_N \subseteq \mathbb X, \ X_k \in \mathcal X, \ X_N \in \mathcal X, \ \mu_k \in \mathcal U} \; .
\end{array} \right.
\end{array} \hspace{-0.4cm}
\end{align}
In the \change{context of this paper, $\mathcal X$ denotes a} suitable class of polytopes---to be specified in Section~\ref{sec::MPCsummary}. This means that the set-valued optimization variables of~\eqref{eq::tocp}, $X_0, X_1, \ldots, X_N$, form a robust forward invariant polytopic tube satisfying the state constraints, $X_k \subseteq \mathbb X$ over the prediction horizon of the controller. In this context, $\hat x$ denotes the current state measurement and the constraint $\hat x \in X_0$ ensures that the first set of our optimized tube contains $\hat x$. \change{Notice that the control laws $\mu = (\mu_0,\mu_1,\ldots,\mu_{N-1})$ are, in 
this formulation, freely optimized. In fact, throughout this article, we shall not impose any restrictions on these control laws apart from requiring that they respect all constraints in~\eqref{eq::tocp}. As such, the maps $\mu_k:\mathbb R^{n_x}\to\mathbb{U}$ could be nonlinear or even discontinuous functions.}

\bigskip
\noindent
A detailed discussion about how to design the stage cost function $L: \mathcal X \times \mathcal U \to \mathbb R$ as well as the initial- and end-cost functions $L_0: \mathcal X \to \mathbb R$ and $L_N: \mathcal X \to \mathbb R$ can be found in Section~\ref{sec::MPCsummary}.

\section{Template Polyhedra}
\label{sec::polyhedra}
Our construction of a tractable reformulation of~\eqref{eq::tocp}
is based on parametric polyhedra, \change{which have been analyzed by many authors, for instance,~\cite{Bemporad2002,Blanchini2008,Mayne2005,Rakovic2012,Rakovic2016}}. Therefore, our introduction of (template) polyhedra in Section~\ref{sec::polyIntro} is kept short, but Sections~\ref{sec::facets} and~\ref{sec::configuration} briefly review how generalized faces and facial configurations of template polyhedra can be defined and analyzed, as such configurations are hardly ever analyzed in the context of optimization and control. At this point, \change{a novel idea is presented in Sections~\ref{sec::confTemplates} and~\ref{sec::vertex}, which introduce a new class of configuration-constrained polytopes that admit a joint parameterization of their facets and vertices. The corresponding main result is summarized in Theorem~\ref{thm::convhull}.}

\subsection{Polyhedra and Polytopes}
\label{sec::polyIntro}
Let $Y \in\mathbb{R}^{m \times n}$ be a given matrix. Throughout this paper, polyhedra with parameter $y \in \mathbb{R}^{m}$ are denoted by
\begin{equation*}
P(Y,y) \, \defeq \, \{ x \in \mathbb{R}^{n} \ | \ Yx \leq y \} \subseteq \mathbb{R}^{n}\;.
\end{equation*}
Every polyhedron is both convex and closed. If $P(Y,y)$ is also bounded, it is called a polytope~\cite{Bitsoris1988,Blanchini2008}.

\begin{definition}
The set of feasible parameters,
\[
\mathbb Y \; \defeq \; \left\{ \ y \in \mathbb R^m \ \middle| \ P(Y,y) \neq \varnothing \ \right\} \; ,
\]
is called the natural domain of the template $Y \in \mathbb R^{m \times n}$.
\end{definition}

\noindent
The set $\mathbb Y$ is unbounded, as the following statement holds independently
of how $Y$ is chosen.

\begin{proposition}
\label{prop::nonempty}
We have $\mathbb R_{+}^m \subseteq \mathbb Y$\change{, where $\mathbb R_{+}^m$ denotes the set of componentwise non-negative vectors in $\mathbb R^m$.}
\end{proposition}

\textbf{Proof.} 
If we have $y \in \mathbb R_+^m$, then the inequality $Yx \leq y$ is trivially satisfied for $x = 0$. Consequently, $P(Y,y)$ is non-empty for all $y \in \mathbb R_{+}^m$, which implies $\mathbb R_{+}^m \subseteq \mathbb Y$.
\qed

\bigskip
\noindent
Apart from the above statement, we will see below that $\mathbb Y$ is itself a polyhedron, as we will show as a side product of the more general considerations in the sections below.

\subsection{Faces of Polyhedra}
\label{sec::facets}
Since the developments in this paper build upon understanding the geometry of polyhedra, this section briefly recalls how faces of polyhedra are defined. Let
\[
I = \{ i_1, i_2, \ldots, i_{|I|} \}\subseteq \{ 1, 2, \ldots m \}
\]
be an index set. In the following, the notation
\begin{align}
Y_I & \defeq \; \left[ \, Y_{i_1}^\tr \, , \, Y_{i_2}^\tr \, , \, \ldots \, , \, Y_{i_{|I|}}^\tr \right]^\tr \notag \\
\text{and} \qquad y_I & \defeq \; \left[ \, y_{i_1} \, , \, y_{i_2} \, , \, \ldots \, , \, y_{i_{|I|}} \right]^\tr \notag
\end{align}
is used to denote the matrix (or vector) constructed by collecting all the rows of $Y$ (or coefficients of $y$), whose index is in the set $I$. Moreover,
\begin{align}
\label{eq::defF}
\mathcal F_I(y) \; \defeq \; \left\{ \ x \in P(Y,y) \ \middle| \ Y_I x \geq y_I \ \right\}
\end{align}
denotes the face of $P(Y,y)$ associated to $I$. This notation is formally also defined for the case $I = \varnothing$, since
\[
\mathcal F_{\varnothing}(y) \; = \;  P(Y,y) \; .
\]
Additionally, we recall that $0$-dimensional faces are called vertices, $1$-dimensional faces are called edges, and faces with co-dimension $1$ are called facets. The intersection of faces is again a face, since
\[
\forall I,J \subseteq \{ 1, 2, \ldots, m\}, \quad \mathcal F_{I}(y) \cap \mathcal F_{J}(y) = \mathcal F_{I \cup J}(y)
\]
Moreover, every face of a polyhedron is a polyhedron. Similarly, the faces of a polytope are also polytopes.

\bigskip
\noindent
Apart from the above standard definitions, we introduce the following notion of entirely simple polyhedra.

\begin{definition}
\label{def::simple}
The polyhedron $P(Y,y)$ is called entirely simple\footnote{We additionally recall that an $n$-dimensional polytope (with $n \geq 1$) is called simple, if all its vertices lie in exactly $n$ edges, or equivalently, if it is dual to a simplicial polytope~\cite{Gruenbaum1967}. This naming convention is consistent with our definition of entirely simple polyhedra in the sense that entirely simple polytopes are---by Definition~\ref{def::simple}---also simple.} if all of its non-empty faces, $\mathcal F_I(y) \neq \varnothing$, satisfy
\[
\rank(Y_I) = |I| \; .
\]
\end{definition}

\noindent
Notice that $n$-dimensional entirely simple polyhedra have the property that all of their non-empty faces $\mathcal F_I(y)$ satisfy $|I| \leq n$, since $n$ is an upper bound on the rank of the template matrix $Y \in \mathbb R^{m \times n}$.

\subsection{Face Configurations}
\label{sec::configuration}
A face configuration of a polyhedron is the collection of index sets of its non-empty faces.

\begin{definition}
The face configuration $\mathcal C(y)$ of $P(Y,y)$ is 
\[
\mathcal C(y) \; \defeq \; \left\{ \ I \subseteq \{ 1, \ldots, m \} \ \middle| \ \mathcal F_I(y) \neq \varnothing \ \right\} \, .
\]
Moreover, a set $\mathcal I \subseteq 2^{\{ 1,2,\ldots, m\}}$ is called a partial face configuration of the polyhedron $P(Y,y)$, if $\mathcal I \subseteq \mathcal C(y)$.
\end{definition}

\begin{figure*}[t]
\centering
\begin{tabular}{ccc}
\begin{minipage}{0.35\textwidth}
\begin{center}
$\mathbf{ P(Y,y)}$ \\[0.2cm]
\includegraphics[width=0.8\textwidth]{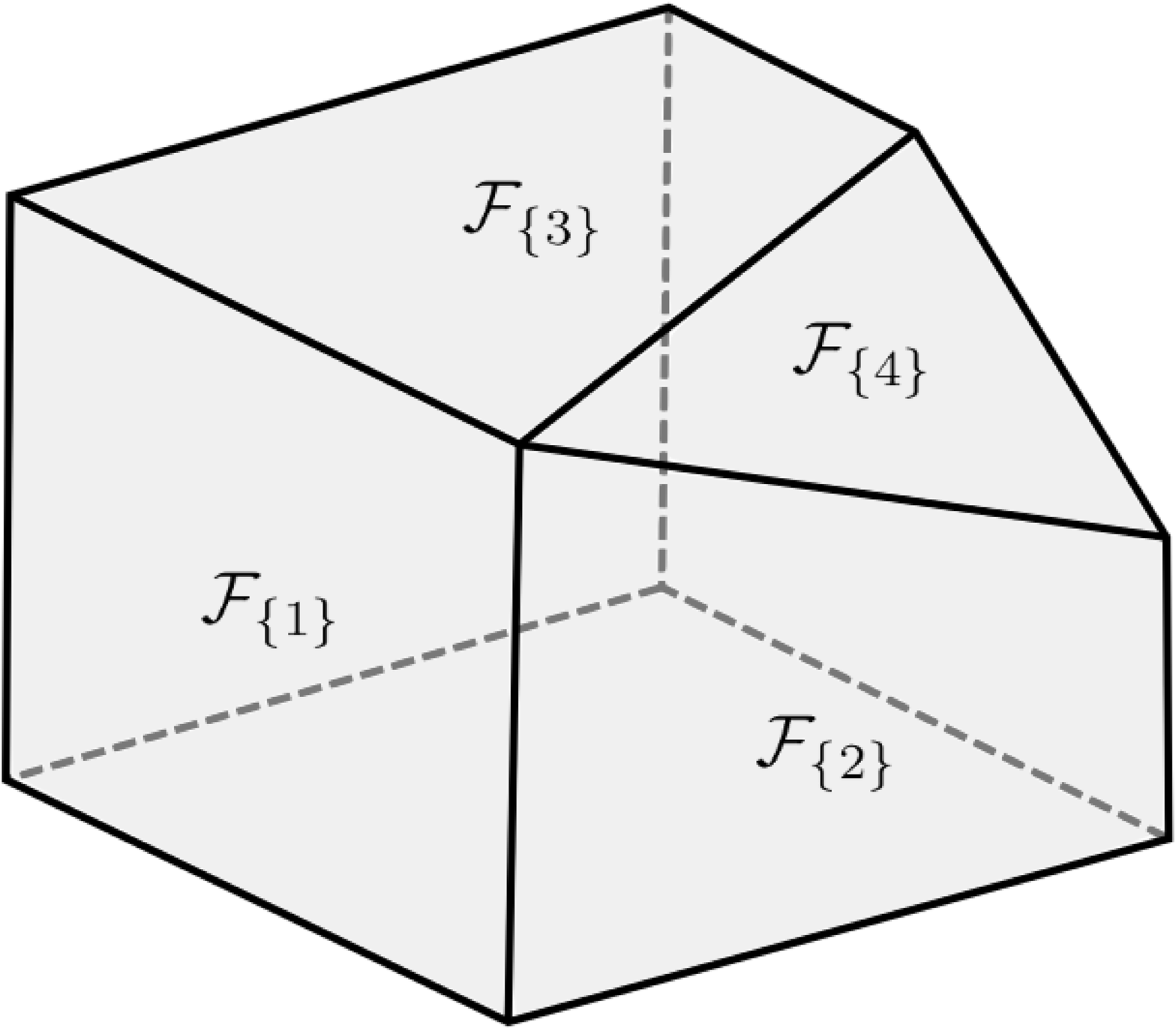}
\end{center}
\end{minipage}
&
\begin{minipage}{0.35\textwidth}
\begin{center}
$\mathbf{ P(Y,y')}$ \\[0.2cm]
\includegraphics[width=0.8\textwidth]{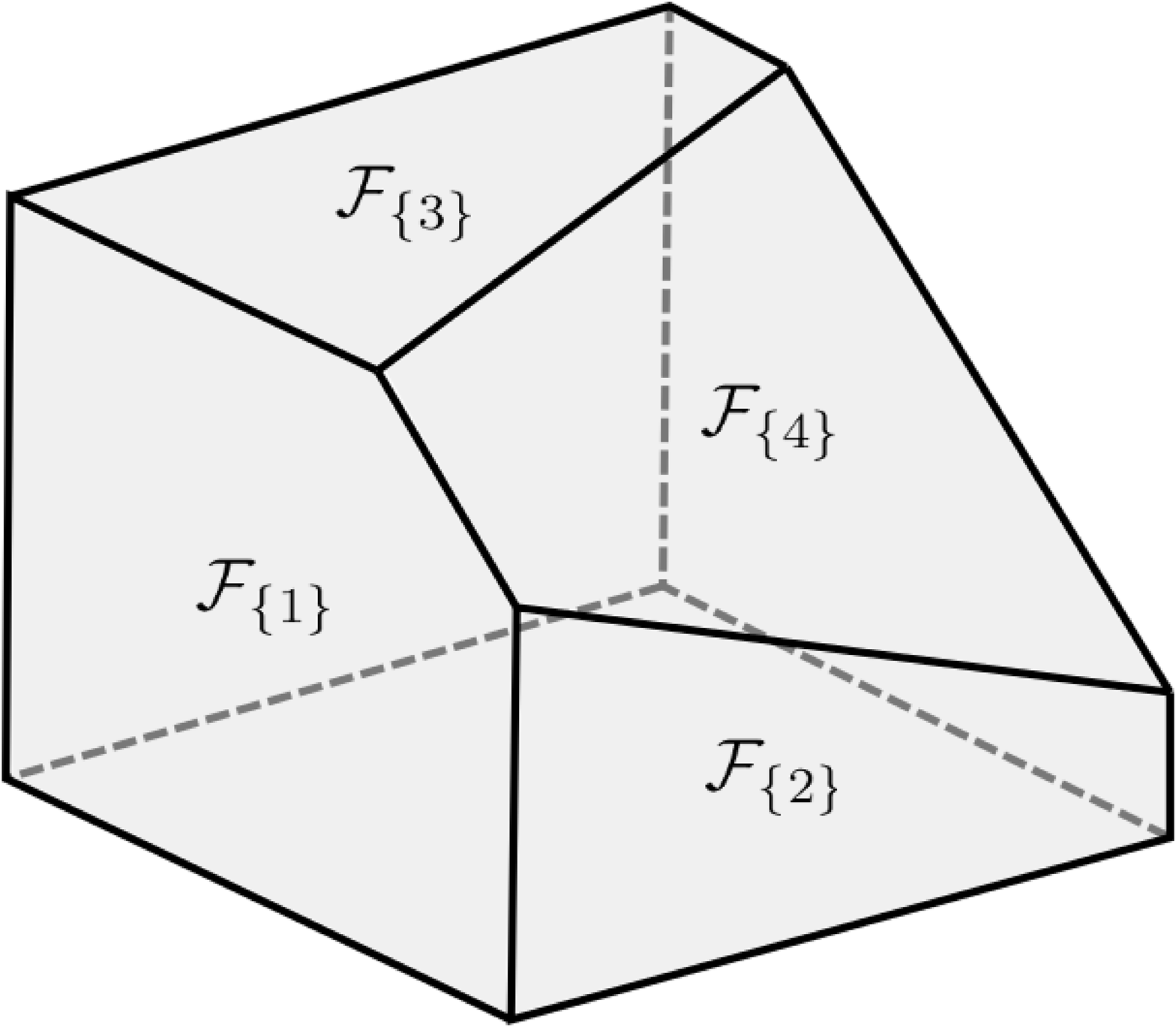}
\end{center}
\end{minipage}
&
\begin{minipage}{0.25\textwidth}
\begin{center}
\textbf{Template matrix}\\[-0.3cm]
\[
Y = \left[
\begin{array}{rrr}
1 & 0 & 0 \\
0 & 1 & 0 \\
0 & 0 & 1 \\
-1 & 2 & 2 \\
-1 & 0 & 0 \\
0 & -1 & 0 \\
0 & 0 & -1 \\
\end{array}
\right]
\]
\vspace{0.2cm}
\end{center}
\end{minipage}
\end{tabular}
\caption{\label{fig::configuration} The polytopes $P(Y,y)$ and $P(Y,y')$ with $y = [1,1,1,3,1,1,1]^\tr$ and $y' = [1,1,1,\nicefrac{5}{2},1,1,1]^\tr$ for a given template matrix $Y \in \mathbb R^{7 \times 3}$. Notice that both polytopes have $7$ two-dimensional facets; that is,
$\mathcal I = \{ \{1 \},\{2 \},\{3 \},\{4 \},\{5 \},\{6 \},\{7\} \}$
is a shared partial face configuration. However, $P(Y,y)$ has $9$ vertices, while $P(Y,y')$ has $10$ vertices. This implies that $P(Y,y)$ and $P(Y,y')$ do not have the same face configuration, $\mathcal C(y) \neq \mathcal C(y')$. 
}
\end{figure*}

\noindent
We say that two polyhedra, $P(Y,y)$ and $P(Y,y')$, have the same configuration if $\mathcal C(y) = \mathcal C(y')$. Similarly, we say that the polyhedra $P(Y,y)$ and $P(Y,y')$ share the partial face configuration $\mathcal I$ whenever
\[
\mathcal C(y) \cap \mathcal I = \mathcal C(y') \cap \mathcal I \; .
\]
At this point, it is recommended to study Figure~\ref{fig::configuration}.
It shows a pair of polytopes whose face configurations do not coincide. The left polytope is not entirely simple, since the face $\mathcal F_I(y)$ with $I = \{ 1,2,3,4 \}$ is non-empty---the corresponding matrix $Y_I \in \mathbb R^{4 \times 3}$ cannot possibly have rank $4$. In contrast to this, the polytope $P(Y,y')$ in the middle of Figure~\ref{fig::configuration} is an entirely simple polytope.

\begin{remark}
\label{rem::PolyConf}
The study of the mathematical properties of polytopes and their face configurations or face lattices has a long history. For instance, the classical formula
\[
\mathsf{v} - \mathsf{e} + \mathsf{f} \; = \; 2 \; ,
\]
relating the number of vertices $\mathsf{v}$, the number of edges $\mathsf{e}$, and number of facets $\mathsf{f}$ of a three dimensional polytope go back to L.~Euler~\cite{Gruenbaum1967}. Similarly, the Dehn-Sommerville relations, which collect similar linear equations that must be satisfied by the numbers of faces of polytopes, have a long history, too~\cite{Bronsted1983}. In general, it turns out to be difficult to count the faces of polytopes. However, proofs of various upper- and lower bound results can be found in the literature from the second half of the $20$th century; see, for example~\cite{Mullen1971,Stanley1980}.
\end{remark}

\bigskip
\noindent
Since there are only finitely many configurations, $\mathcal C$ must be a piecewise constant function on $\mathbb Y$. As such, it is interesting to ask which polyhedra $P(Y,y)$ have a locally stable configuration, such that $\mathcal C(y) = \mathcal C(y')$ holds for all $y'$ in a small open neighborhood of~$y$. As it turns out, such a local stability property holds if and only if $P(Y,y)$ is an entirely simple polyhedron. This statement will be obtained as a special case of a more powerful global stability statement---see Theorem~\ref{thm::regularTemplates} and Corollary~\ref{cor::LocalStability} below.

\subsection{Configuration Domains}
\label{sec::confTemplates}

\change{In} general it is hard to bound the number of faces of high dimensional polyhedra---let alone to compute or to classify them~\cite{Mullen1971,Stanley1980}. \change{Instead}, this paper attempts to characterize the set of polyhedra that share a given (partial) face configuration. In this context, the following definition is useful.

\begin{definition}\label{def::configuration-domain}
Let $\mathcal I$ be a given collection of subsets of the index set $\{ 1,2,\ldots, m\}$. The set
\[
\mathbb Y_{\mathcal I} \defeq \left\{ \ y \in \mathbb R^m \ \middle| \ \forall I \in \mathcal I, \; \mathcal F_I(y) \neq \varnothing \ \right\} \; .
\]
is called the configuration domain of $\mathcal I$.
\end{definition}

\noindent
In words, Definition~\ref{def::configuration-domain} states that a configuration domain is a set of parameters $y$ for which the polyhedron $P(Y,y)$ has certain non-empty faces, specified by the collection $\mathcal I$. An important observation is that configuration domains admit computationally tractable representations. 

\begin{lemma}
\label{lem::cone}
The configuration domain $\mathbb Y_{\mathcal I}$ is for any given set $\mathcal I \subseteq 2^{\{ 1,2,\ldots, m\}}$ a polyhedral cone.
\end{lemma}

\textbf{Proof.}
Let us introduce the auxiliary sets
\begin{align}
\label{eq::YICone}
\mathbb K_I \; \defeq \; \left\{
\ [x^\tr,y^\tr]^\tr \ \middle| \
Y x \leq  y, \; Y_I x \geq y_I \
\right\}
\end{align}
for all $I \in \mathcal I$, which are, by construction, polyhedral cones in $\mathbb R^{n+m}$, since both constraints in~\eqref{eq::YICone} are linear (and, thus, homogeneous) in the stacked vector $[x^\tr,y^\tr]^\tr$. Since the face $\mathcal F_I(y)$ can be written in the form
\[
\mathcal F_I(y) = \{ \ x \in P(Y,y) \ \mid \ Y_I x \geq y_I \ \} \; ,
\]
the set $\mathbb Y_{\{ I \}} = \{ \, y \mid \mathcal F_I(y) \neq \varnothing \, \}$ is the projection of $\mathbb K_I$ onto the last $m$ coordinates; that is, $\mathbb Y_{\{ I \}} = [ \, 0 \; \mathbb{1} \, ] \cdot \mathbb K_I$. Since projections preserve polyhedral  
conic structures, $\mathbb Y_{\{ I \}}$ is a polyhedral cone. Furthermore, since finite 
intersections also preserve polyhedral conic structures, 
\[
\mathbb Y_{\mathcal I} \; = \; \bigcap_{I \in \mathcal I} \mathbb Y_{\{ I \}}
\]
is a polyhedral cone, which completes our proof.
\qed

\bigskip
\noindent
For the special case that the collection $\mathcal I = \{ \varnothing \}$ consists of the empty set only, the configuration domain
\[
\mathbb Y_{\{ \varnothing \}} = \mathbb Y
\]
coincides with the natural parameter domain $\mathbb Y$. Due to the importance for some of the constructions below, we summarize this statement in the following corollary.

\begin{corollary}
\label{cor::nonempty}
The set $\mathbb Y$ is a convex polyhedral cone with non-empty interior in $\mathbb R^m$.
\end{corollary}

\textbf{Proof.}
Since $\mathbb Y = \mathbb Y_{\{ \varnothing \}}$, the fact that $\mathbb Y$ is a polyhedral cone is a special case of Lemma~\ref{lem::cone}. Next, Proposition~\ref{prop::nonempty} states that $\mathbb Y$ contains the open set $\mathbb R_{++}^m$, implying that $\mathbb Y$ has a non-empty interior.
\qed

\bigskip
\noindent
In contrast to $\mathbb Y$, the configuration domains $\mathbb Y_{\mathcal I}$ do not necessarily have a non-empty interior. Therefore, the definition below introduces a notion of regularity for $\mathbb Y_{\mathcal I}$.

\begin{definition}
A configuration domain $\mathbb Y_{\mathcal I}$ is regular if its interior in $\mathbb R^m$ is non-empty.
\end{definition}

\bigskip
\noindent
The theorem below provides a unique characterization of regular configuration domains.

\begin{theorem}
\label{thm::regularTemplates}
Let $\mathcal I$ be a set of subsets of $\{ 1, \ldots, m\}$. The following statements are equivalent:
\begin{enumerate}
\addtolength{\itemsep}{2pt}
\item The set $\mathbb Y_{\mathcal I}$ is a regular configuration domain.
\item There exists a point $\sigma \in \mathbb Y_{\mathcal I}$ such that $P(Y,\sigma)$ is an entirely simple polyhedron.
\end{enumerate}
\end{theorem}

\noindent
Theorem~\ref{thm::regularTemplates} can be interpreted as a global configuration stability result in the sense that it contains the following local stability result as a special case.

\begin{corollary}
\label{cor::LocalStability}
The polyhedron $P(Y,y)$ is locally configuration stable if and only if it is entirely simple.
\end{corollary}

\noindent
Moreover, another important consequence of Theorem~\ref{thm::regularTemplates} is that entirely simple polytopes are prevalent.

\begin{corollary}
\label{cor::Prevalent}
The polyhedron $P(Y,y)$ is entirely simple for almost all $y \in \mathbb Y$.
\end{corollary}

\noindent
Proofs of Theorem~\ref{thm::regularTemplates}, Corollary~\ref{cor::LocalStability}, \change{and Corollary~\ref{cor::Prevalent}} can be found in the Appendices~\ref{sec::proofA} and~\ref{app::B}.

\color{changecolor}

\subsection{Vertex configuration domains of polytopes}
\label{sec::vertex}
Let $P(Y,\sigma)$ be an entirely simple polytope for a given parameter $\sigma \in \mathbb R^m$. We use the notation $\mathcal I_n$ to denote the set of subsets of $\{ 1,2, \ldots, m \}$ with $n$ elements such that
\[
\mathcal V \ \defeq \ \left\{ \ I \in \mathcal I_n \ \middle| \ \mathcal F_I(\sigma) \neq \varnothing \
\right\}
\]
can be interpreted as the index set collection associated with the vertices of $P(Y,\sigma)$. Clearly, if $\overline m$ denotes the number of vertices of $P(Y,\sigma)$, we can enumerate them as
\[
\mathcal V \ = \  \{ \mathcal V_1, \mathcal V_2, \ldots, \mathcal V_{\overline m} \} \quad \text{and define} \quad V_i \ \defeq \ Y_{\mathcal V_i}^{-1} \mathbb{1}_{\mathcal V_i} \; ,
\]
where $\mathbb 1$ denotes the unit matrix and $\mathbb{1}_{\mathcal V_i}$ a matrix with $n$ rows---namely, the $n$ unit vectors of $\mathbb R^m$ whose coordinate index is in $\mathcal V_i$. The matrices $Y_{\mathcal V_i}$ are invertible, since $P(Y,\sigma)$ is assumed to be entirely simple. The points
\[
V_1 \sigma, \ V_2 \sigma, \ \ldots, \ V_{\overline m} \sigma \in P(Y,\sigma)
\]
correspond to the $\overline m$ isolated vertices of $P(Y,\sigma)$. Next, for any given $y \in \mathbb R^m$ an equivalence of the form
\begin{eqnarray}
\mathcal F_{\mathcal V_i}(y) \neq \varnothing &\quad \Longleftrightarrow \quad& V_i y \ \in \ P(Y,y) \notag \\[0.16cm]
&\quad \Longleftrightarrow \quad& Y V_i y \ \leq \ y \notag \\[0.16cm]
\label{eq::FViy}
&\quad \Longleftrightarrow \quad& (Y V_i - \mathbb{1})y \ \leq \ 0
\end{eqnarray}
holds. The latter inequality motivates the introduction of the conic constraint matrix
\begin{align}
\label{eq::E}
E \ \defeq \ \left(
\begin{array}{c}
Y V_1 - \mathbb{1} \\[0.16cm]
Y V_2 - \mathbb{1} \\[0.16cm]
\vdots \\[0.16cm]
Y V_{\overline m} - \mathbb{1}
\end{array}
\right)
\end{align}
such that we can write the vertex configuration domain $\mathbb Y_{\mathcal V}$ in its explicit form
\begin{equation}\label{eq::vertex-configuration-domain}
\mathbb Y_{\mathcal V} \overset{\eqref{eq::FViy},\eqref{eq::E}}{=} \{ \ y \in \mathbb R^{m} \ \mid \ E y \leq 0 \ \} .
\end{equation}
The role of this configuration domain in the ongoing developments of this paper is clarified by the following theorem. It implies that for all $y \in Y_{\mathcal V}$, the polytope $P(Y,y)$ has at most $\overline m$ vertices.

\begin{theorem}
\label{thm::convhull}
Let $P(Y,\sigma)$ be an entirely simple polytope with vertex configuration domain
$\mathbb Y_{\mathcal V}$, as defined in~\eqref{eq::vertex-configuration-domain}. Then, we have
\[
\forall y \in \mathbb Y_{\mathcal V}, \quad  P(Y,y) = \mathrm{conv}( V_1 y, V_2 y, \ldots, V_{\overline m} y ),
\]
where $\mathrm{conv}(\cdot)$ denotes the convex hull operator.
\end{theorem}

\textbf{Proof.} Notice that for the scalar case, $n=1$, we may assume $Y = (-1,1)^\tr$ without loss of generality; that is
$$P[Y,y]= [-y_1,y_2]$$
is an interval. The matrices $V_1 = (-1,0)$ and $V_2 = (0,1)$ locate the two vertices, $V_1 y = -y_1$ and $V_2 y = y_2$. The associated vertex configuration domain is then given by $$\mathbb Y_{\mathcal V} = \{ y \in \mathbb R^2 \mid y_1+y_2 \geq 0 \} \; .$$
Thus, for $n=1$, the statement of the theorem holds, as $P(Y,y) = [-y_1,y_2] = \mathrm{conv}(-y_1,y_2)$ holds for all $y \in \mathbb Y_{\mathcal V}$.

\bigskip
\noindent
Next, we may assume $n \geq 2$. Due to the definition of $\mathbb Y_{\mathcal V}$ the sets $\{ V_i y \}$ are non-empty faces of $P(Y,y)$; that is, $V_i y \in P(Y,y)$ for all $i \in \{ 1, 2, \ldots, \overline m \}$, compare~\eqref{eq::FViy}. Hence, since $P(Y,y)$ is convex, the inclusion
\begin{align}
\label{eq::convInclusion}
\mathrm{conv}( V_1 y, V_2 y, \ldots, V_{\overline m} y ) \subseteq P(Y,y)
\end{align}
holds. In order to establish the reverse inclusion, we need a global geometric argument. Let $\mathcal S^{n-1}$ denote the unit sphere in $\mathbb R^n$. Since we assume $n \geq 2$, the surface area of $\mathcal S^{n-1}$ is well-defined and given by\footnote{We use the notation $\Gamma(n) \defeq \int_0^\infty e^{-t} t^{n-1} \, \mathrm{d}t$ for any $n \geq 1$.}
\[
|\mathcal S^{n-1}| = \frac{2\pi^\frac{n}{2}}{\Gamma(n/2)} \; .
\]
Moreover, for any given point $x \in P(Y,y)$, let
\[
\mathcal N(x,y) \ \defeq \ \left\{ \ z \in \mathcal S^{n-1} \ \middle| \
\begin{array}{l}
\forall x' \in P(Y,y), \\
z^\tr (x'-x) \leq 0
\end{array} \
\right\}
\]
denote the intersection of $\mathcal S^{n-1}$ and the normal cone of $P(Y,y)$ at $x$. 
The Gauss-Federer curvature~\cite{Federer1959} of any vertex $v$ of $P(Y,y)$ 
can now be defined as
\[
\kappa(v,y) \ \defeq \ |\mathcal N(v,y)|,
\]
where $|\mathcal N(v,y)|$ denotes the surface area of the spherical polygon 
$\mathcal N(v,y)$. Next, the key observation of this proof is that $\kappa$ can 
be used as an invariant: if $P(Y,y)$ is an entirely simple polytope with 
$y \in \mathbb Y_{\mathcal V}$, then we have
\begin{align}
\label{eq::alphaInvariance}
\kappa(V_i y,y) = \kappa(V_i \sigma,\sigma) \; .
\end{align}
This follows as the Gauss-Federer curvature of isolated vertices only depends on
the direction of the normal vectors of the facets intersecting at this vertex. Moreover,
\begin{align}
\label{eq::alpha2}
\sum_{i=1}^{\overline m} \kappa(V_i \sigma,\sigma) \ = \ \frac{2\pi^\frac{n}{2}}{\Gamma(n/2)} \; ,
\end{align}
holds because $P(Y,\sigma)$ has exactly $\overline m$ vertices, namely, all points of the form $V_i \sigma$ for $i \in \{ 1,2,\ldots, \overline m \}$. Notice that the above equation can be interpreted as a special case of the Gauss-Bonnet theorem for convex polytopes~\cite{Federer1959}, which, in our case, simply follows from the fact that the spherical polygons $\mathcal N(V_i \sigma,\sigma)$ of the vertices of $P(Y,\sigma)$ form a complete partition of the unit sphere. By substituting~\eqref{eq::alphaInvariance} in \eqref{eq::alpha2}, we find that
\begin{align}
\label{eq::alpha3}
\sum_{i=1}^{\overline m} \kappa(V_i y,y) \ = \ \frac{2\pi^\frac{n}{2}}{\Gamma(n/2)},
\end{align}
still assuming that $P(Y,y)$ is an entirely simple polytope. 
Suppose $P(Y,y)$ had $\overline m + k$ vertices with $k > 0$, where
\[
v_1^\star,v_2^\star,\ldots,v_k^\star \notin \mathrm{conv}( V_1 y, V_2 y, \ldots, V_{\overline m} y )
\]
denote the vertices of $P(Y,y)$ that are not in the convex hull of the $\overline m$ known vertices of the form $V_i y$. If this was the case, then we would have
\[
\sum_{i=1}^{\overline m} \kappa(V_i y,y) + \sum_{i=1}^k \kappa(v_i^\star,y) \ = \ \frac{2\pi^\frac{n}{2}}{\Gamma(n/2)},
\]
which, in turn, would contradict~\eqref{eq::alpha3}, since we have $\kappa(v_i^\star,y) > 0$. Consequently, all vertices of $P(Y,y)$ are contained in the convex hull of the points $V_i y$, \mbox{which---due} \mbox{to~\eqref{eq::convInclusion}---implies} that
\begin{align}
\label{eq::Pconv}
P(Y,y) = \mathrm{conv}( V_1 y, V_2 y, \ldots, V_{\overline m} y ) \; .
\end{align}
Last but not least, due to Corollary~\ref{cor::Prevalent}, the polytope $P(Y,y)$ is entirely simple for almost all $y \in \mathbb Y_{\mathcal V}$. Thus, so far, we have shown that~\eqref{eq::Pconv} holds for almost all $y \in \mathbb Y_{\mathcal V}$. However, since the set of $y$ for which~\eqref{eq::Pconv} holds is closed in $\mathbb R^m$, it follows that the statement of this theorem holds for all $y \in \mathbb Y_{\mathcal V}$.
\qed

\begin{remark}
\label{rem::complexity}
The computational complexity of the above construction depends on the number of facet normals $m$ and the number of vertices $\overline m$ of $P(Y,\sigma)$. For instance, the matrix $E$ has, in general, $\overline m \cdot m$ rows and $m$ columns. There are, however, two important aspects in this context:
\begin{enumerate}
\addtolength{\itemsep}{2pt}

\item The numbers $m$ and $\overline m$ only depend on $Y$ and $\sigma$, which
are, in the constructions below, chosen by us. This choice allows us to trade-off between the accuracy of the set 
representation and its complexity. As an example, hyperboxes in 
$\mathbb R^n$ have $2n$ facet normals and $2^{n}$ vertices. In contrast,
the number of facets and vertices of simplices in $\mathbb{R}^n$ scale linearly 
with $n$.

\item The matrix $E$ is typically sparse. Consequently, its number of non-zero 
entries is much smaller than \mbox{$\overline m \cdot m^2$}. Moreover, $E$ often has a significant number of redundant rows that can be removed using an LP solver; see also the example from Remark~\ref{rem::2D} below.

\end{enumerate}

\end{remark}

\begin{remark}
\label{rem::2D}
In $\mathbb R^2$, one may assume---without loss of generality---that the rows of the matrix $Y$ have the form
\[
Y_i = \left[ \ \cos(\varphi_i), \ \sin(\varphi_i) \ \right] \in \mathbb R^{1 \times 2}
\]
with $0 \leq \varphi_1 < \varphi_2 < \ldots < \varphi_m < 2 \pi$. Notice that one can always normalize the constraints $Y x \leq y$ row-wise and sort them with respect to their argument in polar coordinates. Moreover, we may assume that $m \geq 3$ and
\[
\varphi_{i+1}-\varphi_i < \pi \quad \text{and} \quad \varphi_m - \varphi_1 > \pi
\]
for all $i \in \{ 1, \ldots, m-1\}$, such that $P(Y,y)$ is for all $y \in \mathbb Y$ a polytope. The initial parameter is chosen as $\sigma = [1,1,\ldots,1]^\tr \in \mathbb R^{m}$ and the vertices enumerated as
\[
\mathcal V_i = \{ i, i+1 \} \quad \text{and} \quad \mathcal V_m = \{ m,1\}
\]
for all $i \in \{ 1, \ldots, m-1 \}$ and $\mathcal V = \{ \mathcal V_1, \ldots, \mathcal V_m \}$. Notice that this construction is such that the auxiliary points $\xi_i = Y_i^\tr$ satisfy
\begin{align}
& Y_i \xi_i = \cos(\varphi_i)^2 + \sin(\varphi_j)^2 = 1 = \sigma_i, \notag \\[0.1cm]
\text{and} \quad & Y_j \xi_i = \cos(\varphi_i-\varphi_j) < 1 = \sigma_j \notag
\end{align}
for all $i,j \in \{ 1,\ldots, m\}$ with $i \neq j$. Consequently, we have $\xi_i \in \mathcal F_{\{i\}}$ but $\xi_i \notin \mathcal F_{\{j\}}$ for $i \neq j$. This is sufficient to ensure that $P(Y,\sigma)$ is entirely simple and has \mbox{$\overline m = m$} isolated vertices. It is not difficult to check that its vertex configuration domain is given by
\[
\mathbb Y_{\mathcal V} = \{ y \in \mathbb R^m \mid E y \leq 0 \},
\]
with the conic constraint matrix
\[
E = \left(
\begin{array}{cccccccc}
\Delta_2 & \Sigma_1 & \Delta_1 &  \\
& \Delta_3 & \Sigma_2 & \Delta_2 & \\ 
&  & \ddots & \ddots & \ddots \\ 
& & & \Delta_{m-1} & \Sigma_{m-2} & \Delta_{m-2} \\ 
\Delta_{m-1} & & & & \Delta_m & \Sigma_{m-1} \\
\Sigma_m & \Delta_m & & & & \Delta_1 \\
\end{array}
\right) \in \mathbb R^{m \times m}\;.
\]
All empty spaces in $E$ are filled with zeros while $\Delta_i$ and $\Sigma_i$ denote geometric constants given by
\[
\Delta_i = -\sin(\varphi_{i+1}-\varphi_i) \quad \text{and} \quad \Sigma_i = \sin(\varphi_{i+2}-\varphi_i)
\]
for all $i \in \{ 1, \ldots, m \}$. Here, we have set $\varphi_{m+1} = 2 \pi + \varphi_1$ and $\varphi_{m+2} = 2 \pi + \varphi_2$. This expression for $E$ is constructed using~\eqref{eq::E} and then removing redundant rows. This leads to the above sparse conic constraint matrix, which has $m$ rows and $3m$ non-zero coefficients.\footnote{\change{To be precise, for $n=2$ the number $3m$ is a general upper bound on the number of non-zero coefficients of $E$. There are cases in which $E$ can be further simplified. For instance, for $m=3$ (triangles) it is always possible to find a conic representation with a matrix $E$ that has only one row.}}
We now claim that there exists for every non-empty polytope $P(Y,y')\subset \mathbb R^{2}$ a parameter $y\in\mathbb{R}^{m}$ such that $P(Y,y') = P(Y,y)$ and $Ey\leq 0$. The argument is as
follows. If we set
\begin{equation}\label{eq::ysupport}
\forall i\in\{1,\ldots,m\},
\quad y_{i} \defeq \max_{x\in P(Y,y')} \ Y_{i}x ,
\end{equation}
then $P(Y,y') = P(Y,y)$ holds by construction. Furthermore, due to~\eqref{eq::ysupport}, all facets of $P(Y,y)$ are non-empty. But this is only possible, if the intersections of all neighboring facets are contained in $P(Y,y)$,
\begin{equation*}
\forall i \in \{ 1,\ldots, m \}, \quad V_i y \in P(Y,y)
 \;.
\end{equation*}
This implies that $E y \leq 0$ due to our construction of $E$ and the above claim holds.

\end{remark}

\color{black}

\section{Configuration-Constrained Polytopic Tubes}
\label{sec::Tube}

This section deals with the construction of configuration-constrained polytopic robust forward invariant tubes. \change{As in the previous section, $\mathbb Y_{\mathcal{V}}$ given 
by~\eqref{eq::vertex-configuration-domain}
denotes the vertex configuration domain of a given entirely simple polytope $P(Y,\sigma)$. The matrix $E$ is constructed as in~\eqref{eq::E}. 
Moreover, a sequence of polytopes with parameters $y_0,y_1, \ldots \in \mathbb Y_{\mathcal V}$} is called robust forward invariant if
\begin{align}
\label{eq::yrecursion}
\forall k \in \mathbb N, \qquad P(Y,y_{k+1}) \in F( P(Y,y_k) ) \; .
\end{align}
\change{We recall that $\mathbb U$ is assumed to be closed and convex while $\mathbb W$ denotes a compact set. Throughout the following considerations, we additionally define
\[
\overline w_i \ \defeq \ \max_{w \in \mathbb W} \ Y_i C w
\]
such that $C \mathbb W \subseteq P(Y,\overline w)$ is a tight polytopic enclosure of the compact uncertainty set $C \mathbb W$.}

\subsection{\change{Vertex Control Laws}}
\label{sec::polyTubes}

\color{changecolor}
A well-known necessary and sufficient condition for ensuring that~\eqref{eq::yrecursion} holds is that there exist for every vertex $v_{k,i}$ of the polytope $P(Y,y_k)$ a control input $u_{k,i} \in \mathbb U$ such that the condition
\[
\overline A_j v_{k,i} + \overline B_j u_{k,i} + C w \in P(Y,y_{k+1})
\]
holds for all $w \in \mathbb W$ and all vertices of $\Delta$; that is, for all $j \in \{ 1, \ldots, l\}$. Notice that this statement about the existence of vertex control laws \mbox{has---in} a slightly different \mbox{version---originally} been invented and proven by Gutman and Cwikel~\cite{Gutman1986}. In order to briefly discuss the role of this well-known and historic result in the context of the ongoing developments of the current article, we introduce the convex set
\[
\mathcal F \ \defeq \ \left\{
(y,y^+) \in \mathbb Y_{\mathcal V}^2 \
\middle| \
\begin{array}{l}
\exists u_1,u_2,\ldots,u_{\overline m} \in \mathbb U: \\[0.16cm]
\forall i \in \{ 1, \ldots, \overline m \}, \\[0.16cm]
\forall j \in \{ 1, \ldots, l \}, \\[0.16cm]
Y \overline A_j V_i y + Y \overline B_j u_i + \overline w \leq y^+
\end{array}
\right\}.
\]
For the sake of the completeness, we provide a short proof of the following corollary, although its statement is essentially a direct consequence of Theorem~\ref{thm::convhull}, the original results from~\cite{Gutman1986}, and certain properties of linear systems with multiplicative uncertainties~\cite{Kothare1996}.

\begin{corollary}
\label{cor::F}
The configuration-constrained polytopes $P(Y,y)$ and $P(Y,y^+)$ with $y,y^+ \in \mathbb Y_{\mathcal V}$ satisfy
\begin{align}
\label{eq::RecursivityCond}
P(Y,y^+) \in F(P(Y,y))
\end{align}
if and only if $(y,y^+) \in \mathcal F$.
\end{corollary}

\textbf{Proof.} Our proof is divided into two parts. The first part shows that $(y,y^+)\in \mathcal F$ implies~\eqref{eq::RecursivityCond} and the second part establishes the reverse implication.

\bigskip
\noindent
\textit{Part I.} Let us assume that $(y,y^+)\in \mathcal F$. The definition of $\mathcal F$ implies that $y \in \mathbb Y_{\mathcal V}$ and, consequently, Theorem~\ref{thm::convhull} implies that every $x \in P(Y,y)$ can be written in the form
\[
x = \sum_{i = 1}^{\overline m} \theta_i(x) V_i y
\]
for suitable scalars $\theta_1(x),\theta_2(x),\ldots,\theta_{\overline m}(x) \geq 0$
satisfying $\sum_{i=1}^{\overline m} \theta_i(x) = 1$. Moreover, the definition of $\mathcal F$ implies that there exists for every vertex $V_i y$ a control input $u_i \in \mathbb U$ such that
\begin{align}
\label{eq::vertexLaw}
Y \overline A_j V_i y + Y \overline B_j u_i + \overline w \leq y^+
\end{align}
for all $j \in \{ 1,\ldots, l \}$. The corresponding control law
\begin{align}
\label{eq::VertexControlLaw}
\mu(x) \ \defeq \ \sum_{i=1}^{\overline m} \theta_i(x) u_i
\end{align}
satisfies $\mu(x) \in \mathbb U$, since $\mathbb U$ is convex. Next, for any given $x \in P(Y,y)$, we can multiply~\eqref{eq::vertexLaw} with $\theta_i(x)$ on both sides and take the sum over $i$ in order to show that
\begin{align}
\label{eq::vertexCo2}
\forall j \in \{ 0,\ldots, l \}, \quad Y \left( \overline A_j x + \overline B_j \mu(x) \right) + \overline w \leq y^+ \; .
\end{align}
Moreover, for the same given $x \in P(Y,y)$, we have
\begin{align}
& \forall w \in \mathbb W, \ \forall [A, B] \in \Delta, \quad A x + B \mu(x) + C w \in P(Y,y^+) \notag \\[0.16cm]
& \Longleftrightarrow \quad
\left\{
\begin{array}{l}
\forall i \in \{ 1, \ldots, m \}, \\[0.16cm]
\underset{w \in \mathbb W}{\max} \ \underset{[A, B] \in \Delta}{\max} \ Y_i \left( A x + B \mu(x) + C w \right) \leq y_i^+
\end{array}
\right.
\notag \\[0.16cm]
& \Longleftrightarrow \quad \left\{
\begin{array}{l}
\forall i \in \{ 1, \ldots, m \}, \\[0.16cm]
\underset{[A, B] \in \Delta}{\max} \ Y_i \left( A x + B \mu(x) \right) + \overline w_i \leq y_i^+
\end{array}
\right.
\notag \\[0.16cm]
& \Longleftrightarrow \quad \left\{
\begin{array}{l}
\forall j \in \{ 0,\ldots, l \} \\[0.16cm]
Y \left( \overline A_j x + \overline B_j \mu(x) \right) + \overline w \leq y^+ \; .
\end{array}
\right.
\notag
\end{align}
Finally,~\eqref{eq::vertexCo2} together with the latter equivalence imply that~\eqref{eq::RecursivityCond} holds. 

\bigskip
\noindent
\textit{Part II.} Reversely, if~\eqref{eq::RecursivityCond} holds for given $y,y^+ \in \mathbb Y_{\mathcal V}$, there exists for every point $V_i y \in P(Y,y)$ a control input $u_i \in \mathbb U$ such that for all $[A, B] \in \Delta$ we have
\begin{align}
& \forall w \in \mathbb W, \quad A V_i y + B u_i + C w \in P(Y,y^+) \notag \\[0.16cm]
\Leftrightarrow \quad & \forall j \in \{ 1,\ldots, m \}, \ \ \max_{w \in \mathbb W} \, Y_j(A V_i y + B u_i + C w) \leq y_j^+ \notag \\[0.16cm]
\Leftrightarrow \quad & Y A V_i y + Y B u_i + \overline w \leq y^+ \; . \notag
\end{align}
In particular, the latter equivalence holds at all vertices of the matrix polytope $\Delta$. This is sufficient to conclude that $(y,y^+) \in \mathcal F$.
\qed

\begin{remark}
\label{rem::ControlLaws}
If $P(Y,y)$ is a simplex with $\overline m = n_x+1$ isolated vertices, one can interpolate the vertex control inputs $u_1,u_2,\ldots,u_{\overline m}$ by an affine control law of the form
\[
\mu(x) = K x + k \; .
\]
Here, $K \in \mathbb R^{n_u \times n_x}$ and $k \in \mathbb R^{n_u}$ can be found by solving the linear equation system
\begin{equation}\label{eq::Klinear-interpolation}
\forall i \in \{ 1, \ldots, \overline m \}, \qquad  K V_i y + k = u_i.
\end{equation}
However, if $P(Y,y)$ has $\overline m > n_x +1$ isolated vertices, System~\eqref{eq::Klinear-interpolation} is, in general, overdetermined. Thus, vertex control laws cannot always be interpolated by affine control laws. This observation will be exploited in Section~\ref{sec::caseStudy} to explain why CCTMPC is potentially less conservative than Tube MPC schemes that use affine feedback laws.
\end{remark}

\color{black}

\subsection{Contractivity of Polytopic Tubes}
\label{sec::contractivePolytopes}

\change{Corollary~\ref{cor::F} provides a basis for reformulating the robust 
forward invariance condition~\eqref{eq::yrecursion} as a computationally tractable convex feasibility condition. However, for the development of} a practical Tube MPC formulation, \change{additional assumptions are needed.}

\begin{assumption}
\label{ass::contractive}
The polytope $P(Y,\sigma)$ is entirely simple, \change{feasible,} and $\beta$-contractive, \change{with $\beta < 1$. That is,}
\[
P(Y,\sigma)\subseteq \mathbb X\quad\text{and}\quad \beta \cdot P(Y,\sigma) \in F( P(Y,\sigma) ) \; .
\]
\end{assumption}

\noindent
To begin discussing in which sense Assumption~\ref{ass::contractive} could be considered appropriate for our purposes, one needs to understand first that assuming that practical control systems admit a compact, convex, and $\beta$-contractive set $\mathcal C \subseteq \mathbb X$ for a $\beta < 1$, is not all too restrictive; at least not for the considered class of linear systems with additive and multiplicative uncertainties~\cite[Thm.~7.2]{Blanchini2008}. As soon as one accepts the assumption on the existence of such a set, any entirely simple polytope $P(Y,\sigma)$ with
\begin{align}
\label{eq::POLYapprox}
\sqrt{\beta} \cdot \mathcal C \subseteq P(Y,\sigma) \subseteq \mathcal C
\end{align}
satisfies Assumption~\ref{ass::contractive} after setting $\beta \leftarrow \sqrt{\beta}$. In order to answer questions regarding the worst-case complexity of this construction in dependence on $n_x$ and $\beta$, it is mentioned here that the problem of approximating compact convex sets with polytopes is well-studied in the literature~\cite{Chen2004,Schneider1981}.

\begin{remark}
\label{rem::Y}
Let the polytope $P(Y,\sigma)$ be entirely simple but not necessarily contractive,
and let $\mathcal V$ be its given vertex configuration. Then, there exists a feasible $\beta$-contractive polytope $P(Y,y) \subseteq \mathbb X$ with $y \in \mathbb Y_\mathcal V$ if and only if the convex optimization problem
\begin{align}
\label{eq::feasibilityProblem}
\min_{y,u} \, 0 \quad \mathrm{s.t.} \ \  \left\{
\begin{array}{l}
\forall i \in \{ 0, \ldots, \overline m \}, \ \forall j \in \{ 0, \ldots, l \}, \\[0.16cm]
Y \overline A_j V_i y + Y \overline B_j u_i + \overline w \leq \beta y \\[0.16cm]
E y \leq 0, \ u_i \in \mathbb U, \ V_i y \in \mathbb X
\end{array}
\right.
\end{align}
has a feasible point $(y^\star,u^\star)$. The constraints $V_i y \in \mathbb X$ in~\eqref{eq::feasibilityProblem} enforce feasibility of $P(Y,y)$, since the equivalence
\[
P(Y,y) \subseteq \mathbb X \quad  \Longleftrightarrow \quad 
 \forall i \in \{ 1, \ldots, \overline m \}, \ \ V_i y \in \mathbb X.
\]
holds for all $y \in \mathbb Y_{\mathcal V}$. As such, the above statement is a direct consequence of Corollary~\ref{cor::F}. Thus, if~\eqref{eq::feasibilityProblem} has a feasible solution $(y^\star,u^\star)$, one can simply set $\sigma \leftarrow y^\star$ in order to satisfy Assumption~\ref{ass::contractive}---possibly after adding a small perturbation to $y^\star$ and making $\beta$ slightly larger in the unlikely event that $P(Y,y^\star)$ is not entirely simple. Otherwise, if~\eqref{eq::feasibilityProblem} is infeasible, one
needs to add more rows to $Y$ and repeat the above procedure with a different $\sigma$ until a feasible and strictly contractive polytope is found.
\end{remark}

\noindent
The following lemma is the basis for the construction of stable \change{CCTMPC} controllers.

\begin{lemma}
\label{lem::gamma}
Let Assumption~\ref{ass::contractive} hold, let $y_\mathrm{s} \in \mathbb R^m$ satisfy $y_\mathrm{s} \leq \sigma$ \change{and \mbox{$(y_\mathrm{s},y_\mathrm{s}) \in \mathcal F$}, and let}
\begin{align}
\label{eq::gamma}
\hspace{-0.2cm}
\gamma \; \defeq \; \min_{\gamma' \geq 0} \; \gamma' \quad \mathrm{s.t.} \quad \gamma' (\sigma-y_\mathrm{s}) \geq \beta \sigma - y_\mathrm{s}
\end{align}
denote a contraction constant, $\gamma < 1$. Then, there exists a sequence $y_0,y_1, \ldots \in \mathbb R^m$ such that for all $x_0 \in P(Y,\sigma)$\\[-0.7cm]
\begin{enumerate}
\addtolength{\itemsep}{1pt}
\item the initial feasibility condition, $Y x_0 \leq y_0$, holds,
\item we have $(y_k,y_{k+1}) \in \mathcal F$ for all $k \in \mathbb N$, and
\item the sequence converges $\gamma$-exponentially; that is,
\[
\forall k \in \mathbb N, \qquad \left\| y_k - y_\mathrm{s} \right\|_\infty \leq \gamma^k \, \left\| \sigma - y_\mathrm{s} \right\|_\infty \, . \quad
\]
\end{enumerate}
\end{lemma}

\bigskip
\noindent
\textbf{Proof.} \change{The} key idea of this proof is to show that the particular sequence
\begin{align}
\label{eq::lemClaim}
\forall k \in \mathbb N, \qquad y_k \defeq \gamma^k \sigma + (1 - \gamma^k) y_\mathrm{s}
\end{align}
satisfies all three requirements of the lemma. \change{We} start with the case $k=0$ \change{for which}~\eqref{eq::lemClaim} yields\footnote{We use the definition $0^0 \defeq 1$ for $\gamma = 0$.} $
y_0 = \sigma$. This choice of $y_0$ satisfies the first condition of the lemma, since $x_0 \in P(Y,\sigma) = P(Y,y_0)$ implies $Y x_0 \leq y_0$. Moreover, the definition of $\gamma$ ensures that the inequality
\begin{eqnarray}
y_{k+1} &=& \gamma^{k+1} \sigma + (1-\gamma^{k+1}) y_\mathrm{s} = y_\mathrm{s} + \gamma^k \left[ \, \gamma( \sigma-y_\mathrm{s} ) \, \right] \notag \\
&\overset{\eqref{eq::gamma}}{\geq}&
y_\mathrm{s} + \gamma^k \left[ \, \beta \sigma - y_\mathrm{s} \, \right] =  \gamma^k \beta \sigma + (1-\gamma^k ) y_\mathrm{s} \label{eq::auxyp1}
\end{eqnarray}
holds. \change{In addition, Corollary~\ref{cor::F} and Assumption~\ref{ass::contractive} imply 
$(\sigma,\beta \sigma) \in \mathcal F$ and we recall that $(y_\mathrm{s},y_\mathrm{s}) \in \mathcal F$ holds due to the assumptions of this lemma. Consequently, we have
\begin{align}
\label{eq::FF1}
\gamma^k (\sigma,\beta \sigma) + (1-\gamma^k)(y_\mathrm{s},y_\mathrm{s}) \in \mathcal F,
\end{align}
as $\mathcal F$ is convex. Since $\sigma \in \mathbb Y_{\mathcal V}$ (per Assumption~\ref{ass::contractive}) and $y_\mathrm{s} \in \mathbb Y_{\mathcal V}$ (as $(y_\mathrm{s},y_\mathrm{s}) \in \mathcal F$) holds,~\eqref{eq::lemClaim} implies that
\begin{align}
\label{eq::ykF}
y_{k+1} \in \mathbb Y_{\mathcal V} \; ,
\end{align}
as $\mathbb Y_{\mathcal V}$ is a convex cone. A third consequence of~\eqref{eq::lemClaim} is 
\begin{align}
\label{eq::FF2}
( y_k, \gamma^k \beta \sigma + (1-\gamma^k ) y_\mathrm{s} ) \overset{\eqref{eq::lemClaim}}{=} \gamma^k (\sigma,\beta \sigma) + (1-\gamma^k)(y_\mathrm{s},y_\mathrm{s}).
\end{align}
Thus, after substituting~\eqref{eq::FF2} in~\eqref{eq::FF1} we find that
\[
( y_k, \gamma^k \beta \sigma + (1-\gamma^k ) y_\mathrm{s} ) \in \mathcal F
\]
and then it follows from~\eqref{eq::auxyp1} and~\eqref{eq::ykF} that we have $(y_k,y_{k+1}) \in \mathcal F$.} Finally, the construction of $(y_k)_{k \in \mathbb N}$ in~\eqref{eq::lemClaim} is such that
\[
\left\| y_k - y_\mathrm{s} \right\|_\infty = \left\| \gamma^k ( \sigma - y_\mathrm{s}) \right\|_\infty = \gamma^k \, \left\| \sigma - y_\mathrm{s} \right\|_\infty
\]
implying that the last condition of this lemma is satisfied, too. This completes our proof.
\qed

\section{Polytopic Tube MPC}
\label{sec::MPCsummary}
The goal of this section is to develop a \change{convex} reformulation of~\eqref{eq::tocp}. We specialize on configuration-constrained polytopic tubes with domain
\[
\change{\mathcal X = \{ \ P(Y,y) \ \mid \ y \in \mathbb Y_{\mathcal V} \  \}}
\]
recalling that \change{$\mathbb Y_{\mathcal V}$} denotes \change{the vertex configuration domain of the $\beta$-contractive polytope $P(Y,\sigma)$}.

\subsection{\change{Stage Cost Function}}
\change{Let us assume that the stage cost function has the form}
\begin{align}
\label{eq::Lassumption}
L(P(Y,y),\change{\mu}) &= \; \change{\ell(y,u)} \; ,
\end{align}
where $\ell$ can be any \change{Lipschitz continuous} and strictly convex function in $y$ and \change{$u$. 
Here, we recall that the vertex control law $\mu$ in~\eqref{eq::VertexControlLaw}
depends on the parameter $y$ of the current polytope and the vertex control inputs
$u$. As such, it seems reasonable to assume that $L$ can be expressed as in~\eqref{eq::Lassumption} for a suitable function $\ell$, although the additional assumptions that $\ell$ is strictly convex and Lipschitz continuous are certainly restrictions. These assumptions are, however, motivated by our desire to reformulate~\eqref{eq::tocp} as a convex optimization problem.}

\color{changecolor}

\begin{example}
Let $\mathcal X$ and $V_1, V_2, \ldots, V_{\overline m}$ be defined as before, recalling that the vertices of $P(Y,y)$ are given by $V_i y$ as long as $Ey \leq 0$. The average over the vertices,
\[
\overline y \ \defeq \ \frac{1}{\overline m} \sum_{i=1}^{\overline m} V_i y \ = \ \overline V y \quad \text{with} \quad \overline V \ \defeq \ \frac{1}{\overline m} \sum_{i=1}^{\overline m} V_i,
\]
can be interpreted as the center of $P(Y,y)$. Similarly, if $u_1,u_2,\ldots,u_{\overline m} \in \mathbb U$ denote associated vertex controls,
\[
\overline u \ \defeq \ \frac{1}{\overline m} \sum_{i=1}^{\overline m} u_i \ = \ \overline U u \quad \text{with} \quad \overline U \ \defeq \ \frac{1}{\overline m} \left( \mathbb{1}, \ldots, \mathbb{1} \right)
\]
corresponds to an average control input. An example for a practical quadratic stage cost function is then given by
\begin{eqnarray}
\ell(y,u) &\ = \ & \left\| \overline V y \right\|_{\mathsf{Q}}^2 + \left\| \overline U u \right\|_{\mathsf{R}}^2  \notag \\[0.1cm]
& & + \sum_{i=1}^{\overline m} \left\{ \left\| V_i y - \overline V y \right\|_{\mathsf{S}}^2 + \left\| u_i - \overline U u \right\|_{\mathsf{T}}^2 \right\}, \label{ex::eq::ell}
\end{eqnarray}
where $\mathsf{Q},\mathsf{R},\mathsf{S}$, and $\mathsf{T}$ denote positive definite weighting matrices. Such costs model a trade-off between average tracking performance and the squared distance of the tube vertices and vertex control inputs to their average values.
\end{example}

\color{black}

\subsection{Initial Cost Function}
Formulating Tube MPC controllers which are stable in the enclosure sense requires, in general, an initial cost function \change{$L_0$}, see~\cite{Villanueva2020} for details. Let us construct \change{$L_0$} by introducing the $1$-step cost-to-travel function~\cite{Houska2017}
\begin{eqnarray}
V(y,y^+) &\, \defeq \,& \min_{u} \; \change{\ell(y,u)} \notag \\[0.1cm]
\label{eq::costToTravel}
& & \; \text{s.t.} \; \left\{
\begin{array}{l}
\change{\forall i \in \{ 1,\ldots, \overline m \}, \forall j \in \{ 1, \ldots, l \}, } \\[0.16cm]
\change{ Y \overline A_j V_i y + Y \overline B_j u_i + \overline w \leq y^+,
} \\[0.16cm]
\change{V_i y \in \mathbb X, \ u_i \in \mathbb U, \ V_i y^+ \in \mathbb X}, \\[0.16cm]
\change{E y \leq 0, \ E y^+ \leq 0,}
\end{array}
\right.
\end{eqnarray}
where $V(y,y^+) \defeq \infty$ whenever~\eqref{eq::costToTravel} is infeasible.  Next, we compute an invariant set by solving
\begin{align}
\label{eq::steady}
V_\mathrm{s} \ \defeq \ \min_{y,y^+} \ V(y,y^+) \quad \text{s.t.} \quad 
\change{y = y^+ \, \mid \ \lambda} \; .
\end{align}
Assumption~\ref{ass::contractive} ensures that $(\sigma, \sigma)$ is a feasible point of~\eqref{eq::steady}. Consequently, since $\ell$ is strictly convex,~\eqref{eq::steady} is a feasible \change{and} strictly convex optimization problem. \change{As such, the following assumption is mild\footnote{\change{If Assumption~\ref{ass::contractive} holds, if $\ell$ is a proper, strictly convex, and radially unbounded function, and if Slater's constraint qualification is satisfied, Assumption~\ref{ass::regular} always holds~\cite{BoyVan:ConOpt:04}.}}.
\begin{assumption}
\label{ass::regular}
The strictly convex optimization problem~\eqref{eq::steady} admits a unique primal solution $(y_\mathrm{s},y_\mathrm{s})$ and a dual solution $\lambda \in \mathbb R^{n_x}$, such that strong duality holds.
\end{assumption}
Under Assumption~\ref{ass::regular},} the rotated cost-to-travel function
\begin{align}
\label{eq::RDEF}
R(y,y^+) \ \defeq \ \ V(y,y^+) + \change{\lambda^\tr( y -y^+)}- V_\mathrm{s}
\end{align}
is positive definite with respect to the point $(y_\mathrm{s},y_\mathrm{s})$. Thus, we can define the initial cost function
\[
\forall y \in \mathbb Y, \quad \change{L_0}(P(Y,y)) \; \defeq \; \change{\lambda^\tr y}\;.
\]
This has the advantage that~\eqref{eq::tocp} is equivalent to
\begin{align}
\min_{y} \ \ \ & \sum^{N-1}_{k=0} R(y_k,y_{k+1}) + \change{ \lambda^\tr} y_N   + \change{L_N}(P(Y,y_N)) \label{eq::tocpR} \\[0.16cm]
\notag
\text{s.t.} \ \ \ & \hat x  \in P(Y,y_0) \; ,
\end{align}
which can be interpreted as a positive definite tracking problem in the parameter space. To achieve stability in the enclosure sense, we must ensure that the \change{parameter sequence $y_k$ converges to $y_\mathrm{s}$}---see~\cite{Villanueva2020} for details.

\subsection{Terminal Cost Function}
In this section, a stabilizing terminal cost is constructed. \change{For this aim}, we assume that \change{the initial polytope $P(Y,\sigma)$ is large enough to contain the optimal invariant set,}
\begin{align}
\label{eq::terminalSet}
y_\mathrm{s} \leq \sigma \; .
\end{align}
Moreover, \change{we solve~\eqref{eq::costToTravel} offline at $(y,y^+) = (\sigma,\sigma)$ and at $(y,y^+) = (y_\mathrm{s},y_\mathrm{s})$. The corresponding minimizers are denoted by $u_\sigma$ and $u_\mathrm{s}$, such that
\begin{eqnarray}
\ell(\sigma,u_{\sigma}) &\ =\ & V(\sigma,\sigma) \ = \ R(\sigma,\sigma) + V_\mathrm{s} \\[0.16cm]
\text{and} \qquad \ell(y_\mathrm{s},u_{\mathrm{s}}) &\ =\ & V(y_\mathrm{s},y_\mathrm{s}) \ = \  V_\mathrm{s} \; .
\end{eqnarray}
Let $\overline \ell$ denote the Lipschitz constant of $\ell$ on $P(Y,\sigma)$. In the following, we introduce the shorthand
\begin{align}
\label{eq::rho}
\rho \ \defeq \ \overline \ell \left\| \left(
\begin{array}{c}
\sigma - y_\mathrm{s} \\
u_\sigma - u_\mathrm{s}
\end{array}
\right) \right\|_2 + (1-\gamma) \left| \lambda^\tr(\sigma - y_{\mathrm{s}}) \right|,
\end{align}
which can be pre-computed offline. Here, $\gamma$ is computed as in Lemma~\ref{lem::gamma}, see Equation~\eqref{eq::gamma}.} Next, our goal is to show that the cost function
$$\change{L_N}(P(Y,y)) \; \defeq \; -\change{\lambda}^\tr y + M(y)$$
with auxiliary function
\begin{eqnarray}
\label{eq::M}
M(y) \, \defeq \, & \min_{\change{\alpha \in [0,1]}}& \; \change{R(y,y_\mathrm{s}+\alpha(\sigma - y_\mathrm{s})) +  \frac{\rho \cdot \alpha}{1 - \gamma}} 
\end{eqnarray}
can be used as stabilizing terminal cost. \change{Notice that $M$ is convex, since $R$ is convex in both arguments.} The following theorem establishes that $M$ is a Lyapunov function.

\begin{theorem}
\label{thm::Lyapunov}
Let Assumptions~\ref{ass::contractive} and~\ref{ass::regular} be satisfied and let $M$ be defined as in~\eqref{eq::M}. The function $M$ is non-negative and satisfies $\change{M(y)} = 0$ if and only if $y = y_\mathrm{s}$. Moreover, the Lyapunov descent condition
\[
\change{\min_{y^+} \, \{  R(y,y^+) + M(y^+) \} \; \leq \; M(y)}
\]
holds for all \change{$y \in \mathbb R^m$.} 
\end{theorem}

\textbf{Proof.}
\change{Notice that~\eqref{eq::M} implies $M(y) \geq 0$ for all $y \in \mathbb R^m$, since $R$ is positive definite and $\frac{\rho \cdot a}{1-\gamma}$ is non-negative. Moreover, $M(y_\mathrm{s}) = 0$. This follows} from~\eqref{eq::M} by substituting $\alpha = 0$ and using \change{that $R(y_\mathrm{s},y_\mathrm{s})= 0 $}. Our next goal is to show that $M(y) > 0$ for all $y \neq y_\mathrm{s}$. In order to see this, assume that $M(y) = 0$. In this case, we must have $\alpha = 0$, because, otherwise, the second term in~\eqref{eq::M} would be strictly positive. Consequently, we must have
$
R(y,y_\mathrm{s}) = 0,
$
but this is impossible for $y \neq y_\mathrm{s}$, since $y_\mathrm{s}$ is---due to Assumption~\ref{ass::regular}---the unique minimizer of $\min_{y} R(y,y_\mathrm{s})$. Thus, $M$ is a positive definite function.

\bigskip
\noindent
Next, in analogy to the construction in the proof of Lemma~\ref{lem::gamma}, for an arbitrary parameter $\alpha \in [0,1]$, we introduce the auxiliary sequences
\begin{align}
\label{eq::ytilde}
\widetilde y_k \; & \defeq \; y_\mathrm{s} + \alpha \gamma^k (\sigma - y_\mathrm{s}) \\[0.16cm]
\label{eq::utilde}
\change{\text{and} \qquad \widetilde u_k} \; & \change{\defeq \; u_\mathrm{s} + \alpha \gamma^k ( u_\sigma - u_\mathrm{s}) \; ,}
\end{align}
\change{such that $(\widetilde y_k, \widetilde y_{k+1}) \in \mathcal F$. Since $\widetilde u_k$ is a feasible point of~\eqref{eq::costToTravel} at $(y,y^+) = (\widetilde y_k, \widetilde y_{k+1})$, we find that
\begin{align}
\label{eq::Vbound}
V(\widetilde y_k, \widetilde y_{k+1}) \ \leq \ \ell(\widetilde y_k, \widetilde u_k) \; .
\end{align}
This inequality leads us to the upper bound
\begin{eqnarray}
R(\widetilde y_0, \widetilde y_{1}) &\overset{\eqref{eq::RDEF}}{=}& V(\widetilde y_0, \widetilde y_{1}) + \lambda^\tr( \widetilde y_0 - \widetilde y_{1} ) - V_\mathrm{s} \notag \\[0.16cm]
&\overset{\eqref{eq::Vbound}}{\leq}& \ell( \widetilde y_0, \widetilde u_0) + \lambda^\tr( \widetilde y_0 - \widetilde y_{1} ) - \ell(y_\mathrm{s},u_\mathrm{s}) \notag \\[0.2cm]
&\leq& \overline \ell \cdot  \left\| \left(
\begin{array}{c}
\widetilde y_0 - y_\mathrm{s} \\
\widetilde u_0 - u_\mathrm{s}
\end{array}
\right) \right\|_2 + \left| \lambda^\tr( \widetilde y_0 - \widetilde y_{1}) \right| \notag \\
&\hspace{-0.1cm} \overset{\eqref{eq::ytilde},\eqref{eq::utilde}}{\leq} \hspace{-0.1cm}&
\alpha \cdot \overline \ell \cdot \left\| \left(
\begin{array}{c}
\sigma - y_\mathrm{s} \\
u_\sigma - u_\mathrm{s}
\end{array}
\right) \right\|_2 \notag \\
& & + \, \alpha \cdot \left| \lambda^\tr(\sigma - y_\mathrm{s} ) - \gamma \lambda^\tr(\sigma - y_\mathrm{s} )  \right|
\notag \\
\label{eq::RRbound}
&\overset{\eqref{eq::rho}}{=}& \alpha \cdot \rho \; , 
\end{eqnarray}
which, in turn, implies that the inequality 
\begin{eqnarray}
M(\widetilde y_0) &\ \overset{\eqref{eq::M}}{=} \ & \min_{\alpha' \in [0,1]} \, \left\{ R(\widetilde y_0, y_\mathrm{s} + \alpha'(\sigma - y_\mathrm{s})) + \frac{\rho \cdot \alpha'}{1-\gamma} \right\} \notag \\[0.16cm]
&\leq& R(\widetilde y_0, y_\mathrm{s} + \alpha \gamma (\sigma - y_\mathrm{s})) + \frac{\rho \cdot \alpha \cdot \gamma}{1-\gamma} \notag \\[0.16cm]
&\ \overset{\eqref{eq::ytilde}}{=} \ & R(\widetilde y_0, \widetilde y_1) + \frac{\rho \cdot \alpha \cdot \gamma}{1-\gamma} \ \overset{\eqref{eq::RRbound}}{\leq} \ \frac{\rho \cdot \alpha}{1-\gamma}
\label{eq::My0bound}
\end{eqnarray}
holds for any choice of $\alpha \in [0,1]$ (keeping in mind that $\widetilde y_0$ depends on that value of $\alpha$, too). Consequently, if we choose $\alpha$ as in~\eqref{eq::M} we find that 
\begin{eqnarray}
M(y) & \ \overset{\eqref{eq::M}}{=} \ & R(y,\widetilde y_0) + \frac{\rho \cdot \alpha }{1-\gamma} \notag \\[0.16cm]
& \ \overset{\eqref{eq::My0bound}}{\geq} \ & R(y,\widetilde y_0) + M(\widetilde y_0) \notag \\[0.16cm]
& \ \geq \ & \min_{y^+} \; \{ \, R(y,y^+) + M(y^+) \} \; .
\end{eqnarray}
T}he latter inequality corresponds to the Lyapunov descent statement of the theorem. \qed

\color{changecolor}

\begin{remark}
Notice that existing attempts to construct terminal cost functions for HTMPC~\cite{Rakovic2012} or ETMPC~\cite{Rakovic2016} rely on the availability of robustly stabilizing affine feedback laws. This is in contrast to the above construction of $M$, which does not rely on the explicit availability of such an affine feedback law. In fact, for the proposed construction of $M$, it is not even required that a robustly stabilizing affine feedback law exists in the first place.
\end{remark}

\begin{remark}
The Lyapunov function $M$ is generic in the sense that no assumptions on $\ell$ apart from strict convexity and Lipschitz continuity are needed. However, the bound of $R$ in~\eqref{eq::RRbound} is not sharp. As such, $M$ is not necessarily the best possible approximation of the infinite horizon control performance that could, however, be derived by using sharper cost estimates. For instance, if $\ell$ is known to be a strongly quadratic form, it is not difficult to derive bounds that are tighter than~\eqref{eq::RRbound}.
\end{remark}

\color{black}

\subsection{Implementation of \change{CCTMPC}}
\label{sec::implementation}
\change{Since $\mathbb X$, $\mathbb U$, $\ell$ and $M$ are convex,~\eqref{eq::tocp} is equivalent to the convex optimization problem}
\begin{eqnarray}
&\min_{y,u} \ & \; \lambda^\tr y_0 + \sum^{N-1}_{k=0} \ell(y_k,\change{u_k}) - \lambda^\tr y_N + \change{M(y_N)} \notag \\[0.16cm]
\label{eq::qp}
&\text{s.t.} \ & \; \left\{ \; 
\begin{array}{l}
\forall k \in \{ 0, \ldots, N-1\}, \\[0.16cm]
\change{\forall i \in \{ 1, \ldots, \overline m \}, \ \forall j \in \{ 1, \ldots, l \}, } \\[0.16cm]
\change{Y \overline A_j V_i y_k + Y \overline B_j u_{k,i} + \overline w \leq y_{k+1},} \\[0.16cm]
\change{E y_k \leq 0, \ u_{k,i} \in \mathbb U}, \ \change{V_i y_k \in \mathbb X,} \ \change{Y \hat x \leq y_0} \ .
\end{array}
\right.
\end{eqnarray}
\change{Notice that the} optimization variables of~\eqref{eq::qp} are the \change{tube parameters $y = [y_0,y_1,\ldots,y_N]$ and the vertex control inputs $u = [u_0,u_1,\ldots,u_{N-1}]$.}

\color{changecolor}

\begin{remark}\label{rem::complexityQP}
Assuming that $\mathbb U$ and $\mathbb X$ are polyhedra with $n_{\mathbb U}$ and
$n_{\mathbb X}$ facets respectively, Problem~\eqref{eq::qp} has
\[
(N+1) \cdot m + N \cdot \overline m \cdot n_u 
\]
optimization variables as well as
\[
m \cdot \overline m \cdot l \cdot N + m_E \cdot N + \overline m \cdot (n_{\mathbb X} + n_{\mathbb U}) \cdot N + m
\]
inequality constraints\footnote{\change{In order to be precise, it should be mentioned that, depending on how the function $M$ is represented, additional auxiliary optimization variables and constraints are needed to arrive at a practical implementation; see~\eqref{eq::M} and~\eqref{eq::costToTravel}.}}, where $m_E \ll m \cdot \overline m$ denotes the number of
non-redundant rows of $E$; see Remark~\ref{rem::complexity} for details. As such, under the mentioned assumptions, one may state that~\eqref{eq::qp} has polynomial
run-time complexity. However, it should also be clear that solving~\eqref{eq::qp} is computationally demanding, since the product of $m$, $\overline m$, $l$, and $N$ is potentially a very large number.
\end{remark}

\color{black}

\subsection{Stability Analysis}
Let \change{$y^\star(\hat x)$ and $u^\star(\hat x)$} be parametric minimizers
of~\eqref{eq::qp} in dependence on the measurement $\hat x$. \change{Next, due to the construction of $y^\star$ and $u^\star$, any control law of the form
\begin{align}
\mu_{\mathrm{MPC}}(\hat x) \ \defeq \ \sum_{i=1}^{\overline m} \theta_i(\hat x) u_{0,i}^\star(\hat x)
\end{align}
can be used as a feasible MPC feedback control law, as long as the scalar coefficient functions $\theta_i$ satisfy
\begin{align}
\hat x = \sum_{i=1}^{\overline m} \theta_i(\hat x) V_i y_0^\star (\hat x), \ \ \theta(\hat x) \geq 0 \ \ \text{and} \ \ \left\| \theta(\hat x) \right\|_{1} = 1 \, , \notag
\end{align}
recalling the construction of $\mu$ in the proof of Corollary~\ref{cor::F}, see Equation~\eqref{eq::VertexControlLaw}. For instance, if one wishes to pick a control law with minimal Euclidean norm, one can find suitable coefficients by solving the convex QP
\begin{eqnarray}
\theta(\hat x) \ \in \ &\underset{\theta}{\mathrm{argmin}}& \ \left\| \sum_{i=1}^{\overline m} \theta_i u_{0,i}^\star(\hat x) \right\|_2^2 \notag \\
&\mathrm{s.t.}& \ \left\{ 
\begin{array}{l}
\hat x = \sum_{i=1}^{\overline m} \theta_i V_i y_0^\star (\hat x), \\
\theta \geq 0, \ \sum_{i=1}^{\overline m} \theta_i = 1 \; ,
\end{array}
\right.
\end{eqnarray}
but using other convex control penalty functions is, of course, possible, too.} The associated uncertain closed-loop system of~\eqref{eq::qp} has the form
\[
\forall k \in \mathbb N, \qquad \change{x}_{k+1} = A\change{_k} \change{x}_k + B\change{_k} \mu_{\mathrm{MPC}}(\change{x}_k) + C w_k \; ,
\]
which is started at a given initial state $\change{x}_0 \in \mathbb R^{n_x}$. \change{The following
theorem establishes recursive feasibility and asymptotic stability, in the enclosure
sense, of the controller---see \cite{Villanueva2020} for an in depth discussion 
of set-theoretic stability in Tube MPC schemes.}

\begin{theorem}
\label{thm::stability}
Let Assumptions~\ref{ass::contractive} and~\ref{ass::regular} be satisfied, let \change{$(y_\mathrm{s},\lambda)$} be given by solving~\eqref{eq::steady}, and let \change{$M$} be given by~\eqref{eq::M}. Then the following statements hold independently of the sequences $w_0,w_1, \ldots \in \mathbb W$ \change{and $[A_0,B_0],[A_1,B_1],\ldots \in \Delta$.}
\begin{enumerate}
\addtolength{\itemsep}{2pt}

\item If~\eqref{eq::qp} is feasible for the initial state, $\hat x = \change{x}_0$, \change{then it is also} feasible for all future measurements, $\hat x = \change{x}_k$.

\item We have $\change{x}_k \in P(Y, y_0^\star(\change{x}_k) )$ for all $k \in \mathbb N$.

\item The \change{sequence $y^\star(\change{x}_k)$ is} asymptotically stable,
\begin{eqnarray}
\lim_{k \to \infty} \; y^\star(\change{x}_k) &=& [y_\mathrm{s},y_\mathrm{s}, \ldots, y_\mathrm{s}] \; . \notag
\end{eqnarray}
\end{enumerate}
\end{theorem}

\bigskip
\noindent
\textbf{Proof.} Let us introduce the non-negative function
\[
\mathcal L(\change{y}) \; \defeq \; \sum_{i=0}^{N-1} \change{R(y_i,y_{i+1}) +  M(y_N)},
\]
where $R$ and $M$ are defined as in~\eqref{eq::RDEF} and~\eqref{eq::M}. Moreover, we have $\mathcal L(y) = 0$ if and only if $y = y_\mathrm{s}$, which follows from Theorem~\ref{thm::Lyapunov} and the positive definiteness of $R$. Now, let us assume that~\eqref{eq::qp} is feasible at $\hat x = \change{x}_k$ such that $y^\star(\change{x}_{k})$ \change{is} well-defined. In this case, we can construct the shifted variable \change{sequence}
\begin{align}
\tilde y \; &\defeq [ y_1^\star(\change{x}_k), \ldots, y_N^\star(\change{x}_{k}), y_N^+ ] \; ,
\end{align}
where the last element $y_N^+$ of \change{$\tilde y$} is chosen such that
\begin{align}
\label{eq::LDAUX}
\change{R(y_N^\star(\change{x}_{k}),y_N^+) + M(y_N^+) \leq M(y_N^\star(\change{x}_{k}) )} \; .
\end{align}
The existence of such a vector $y_N^+$ is ensured by Theorem~\ref{thm::Lyapunov}. Moreover, \change{a} feasible shifted sequence for the \change{vertex control variables} of~\eqref{eq::qp} \change{is given by}
\change{
\begin{align}
\tilde u \; &\defeq [ u_1^\star(\change{x}_k), \ldots, u_{N-1}^\star(\change{x}_k), \hat u^\star ] \; ,
\end{align}
where} \change{$\hat u^\star$ is the optimal vertex control inputs that needs to be computed} when evaluating $R(y_N^\star(\change{x}_{k}),y_N^+)$ (see~\eqref{eq::costToTravel} and~\eqref{eq::RDEF}). This construction is such that the shifted variables \change{$(\tilde y, \tilde u)$} are feasible for~\eqref{eq::qp} at the next time instance, because $\hat{x} = x_{k+1}$ satisfies
\[
\change{x}_{k+1} \in P(Y, y_1^\star(\change{x}_k) )\;,
\]
independently of $w_k$ (due to the above definition of $\mu_{\mathrm{MPC}}$), and hence $Y \hat{x} \leq \tilde{y}_0 = y_1^\star(x_k)$. That is, $\hat{x}$ is contained in the robust forward invariant tube. Thus, Statement~1) holds independent of the uncertainty realization, \change{and, additionally, it follows that Statement~2) holds. Moreover}, we have
\begin{eqnarray}
\label{eq::Ld}
\mathcal L( \change{y^\star(\change{x}_k)} ) &\; \overset{\eqref{eq::LDAUX}}{\geq} \;& r_k + \mathcal L( \change{\tilde y} ) \; \geq \; r_k + \mathcal L( \change{y^\star(\change{x}_{k+1})} ) \; ,
\end{eqnarray}
where we have introduced the shorthand
\[
r_k \ \defeq \ R( \change{y_0^\star(\change{x}_k), y_1^\star(\change{x}_k)} ) \; \geq 0 \; .
\]
\change{The last inequality in~\eqref{eq::Ld} follows from the optimality of $y^\star(x_{k+1})$, since our definition of $\mathcal{L}$ is such that this function coincides with the optimal value function of our MPC controller.}
Because $r_k$ is non-negative and because we have \mbox{$r_k = 0$} if and only if $y_0^\star(\change{x}_k) = y_\mathrm{s}$, the function $\mathcal L$ in~\eqref{eq::Ld} is a strictly descending Lyapunov function.  This implies that Statement~3) holds, which completes
the proof. \qed

\section{\change{Tutorial Examples and Numerical Illustration}}
\label{sec::caseStudy}

This section \change{compares CCTMPC with RTMPC, HTMPC, ETMPC, DAFMPC, and FTPMPC. 
Before we present tutorial examples and numerical illustrations, our main observations are summarized.}

\color{changecolor}

\begin{itemize}

\addtolength{\itemsep}{2pt}

\item CCTMPC is never more conservative than polytopic RTMPC: one can always use the invariant cross-section of the rigid tube to construct a vertex configuration domain for CCTMPC. This implies that CCTMPC admits a more flexible tube representation.

\item Using the same argument, CCTMPC is never more conservative than HTMPC: the vertex configuration constraint is invariant under homothetic scaling, since $\mathbb Y_{\mathcal V}$ is a convex cone.

\item For $n_x = 2$, CCTMPC is never more conservative than ETMPC as long as the vertex configuration domain is constructed as in Remark~\ref{rem::2D}. In general, for $n_x \geq 3$, however, the constraint $E y \leq 0$ leads to an actual restriction on the class of representable polytopes. This is in contrast to ETMPC, which does not introduce such a configuration constraint. On the other hand, the construction of ETMPC controllers relies on an affine parameterization of the feedback control policy, while CCTMPC introduces no such restriction; see Remark~\ref{rem::ControlLaws}.
Moreover, existing set-propagation approaches for ETMPC are conservative, complicating the comparison even further. Nevertheless, numerical comparisons indicate that CCTMPC is less conservative than ETMPC---at least for the benchmark examples in $\mathbb R^3$ and $\mathbb R^4$ that we will propose below.

\item It is possible to construct systems with four states for which CCTMPC is systematically less conservative than any separable state feedback controller such as DAFMPC and FPTMPC; see Section~\ref{sec::ex3}.

\end{itemize}

\bigskip
\noindent
In addition to above comments about conservatism, it should be pointed out that RTMPC, HTMPC, ETMPC, and CCTMPC have in common that their computational complexity is given by $\mathcal O(N)$, recalling that $N$ denotes the prediction horizon. In contrast to this, the computational complexity of DAFMPC and FPTMPC is given by $\mathcal O(N^2)$, due to an entirely different parameterization strategy.

\subsection{\change{Example 1}}
\label{sec::ex1}

\begin{figure*}[t]
\begin{center}
\includegraphics[width=0.47\textwidth]{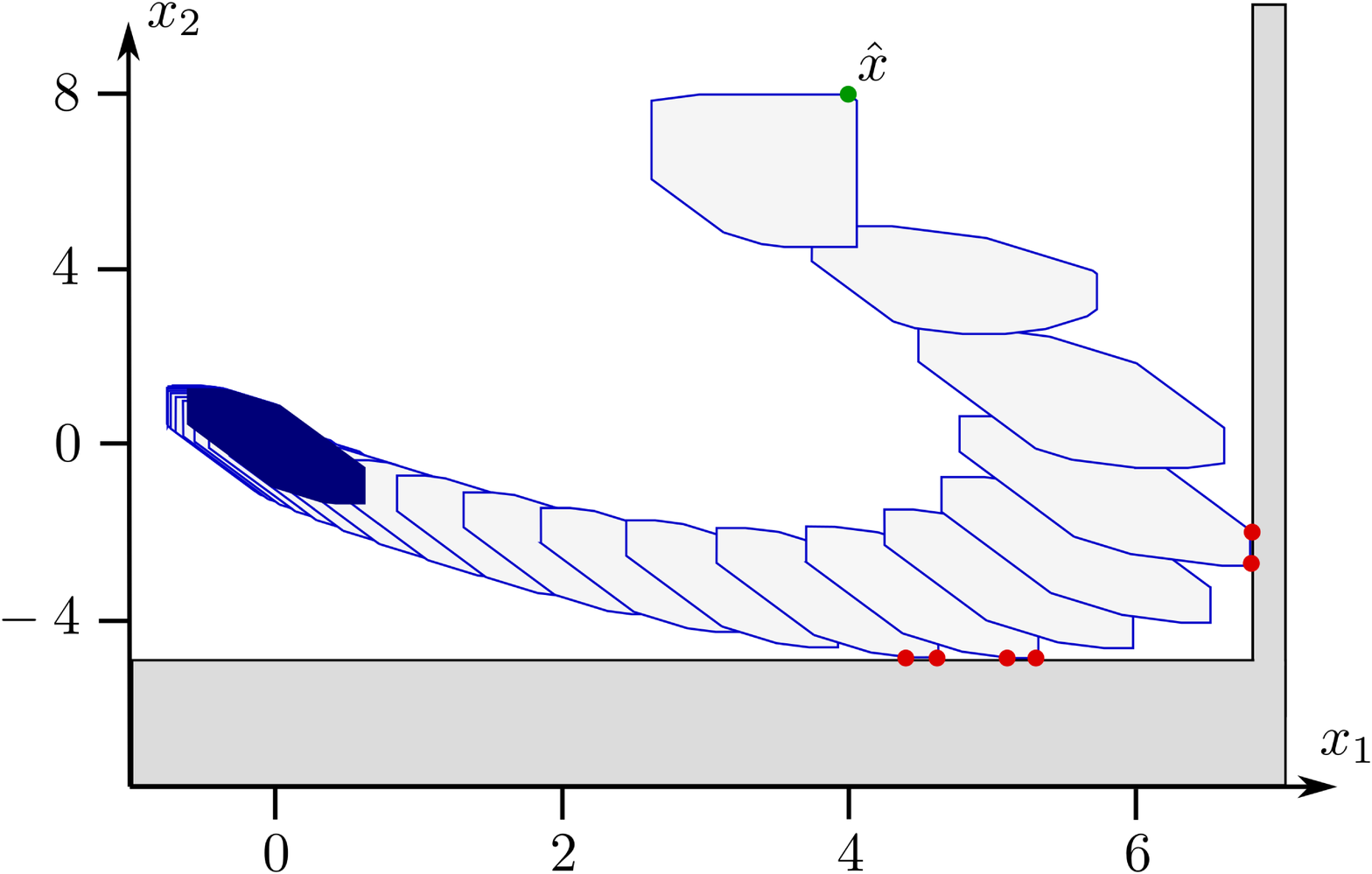}
\hspace{1em}
\includegraphics[width=0.47\textwidth]{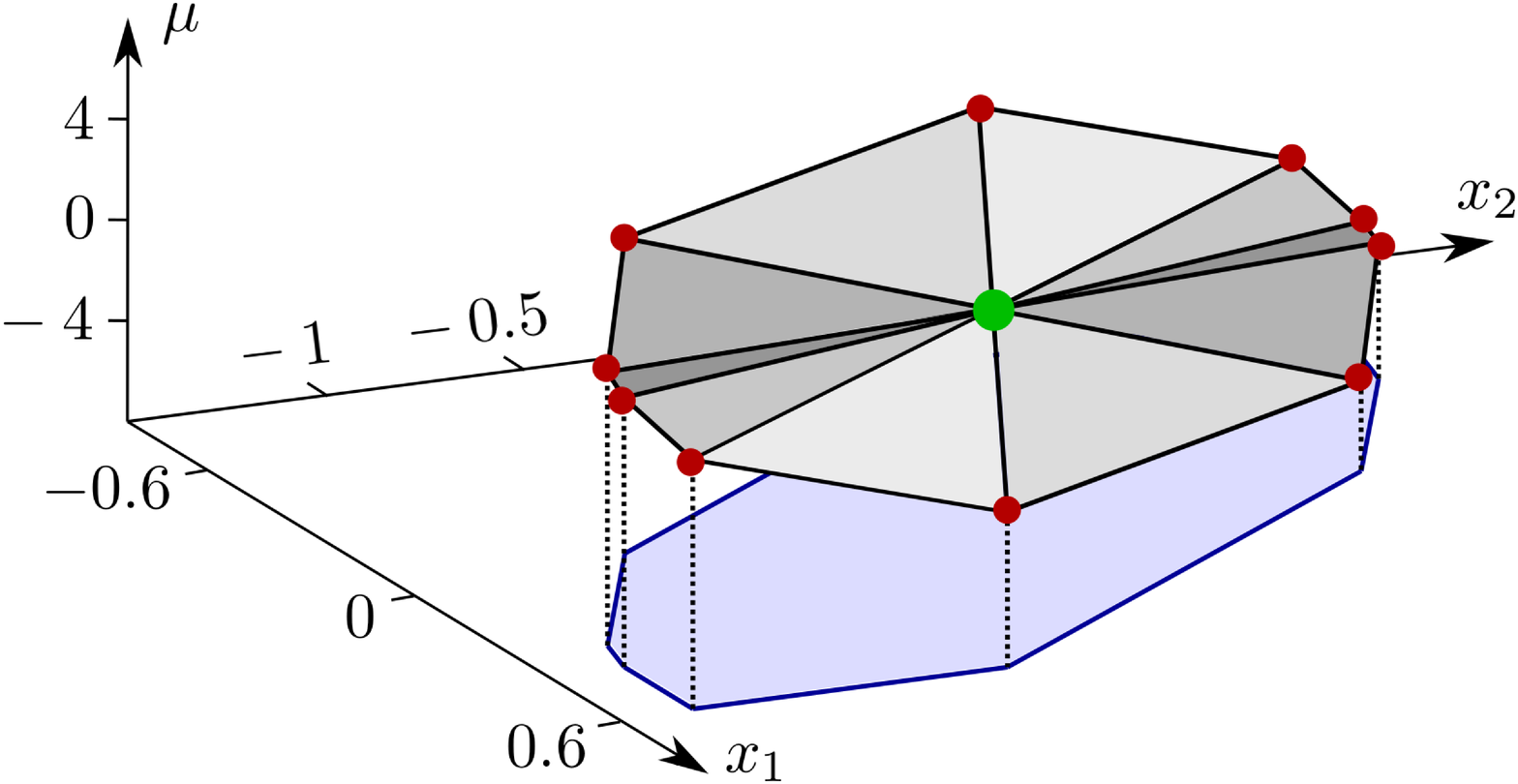}
\end{center}

\caption{\label{fig::prediction} LEFT: Predicted tube (light gray shaded polytopes with blue boundary) \change{for $\hat x = [4,8]^\tr$ (green dot) and $m = 16$}. \change{The} prediction horizon has been set to $N = \change{50}$. The optimal invariant set $P(Y,y_\mathrm{s})$ is indicated by the dark blue shaded set. \change{The} red dots at the \change{fourth,} seventh\change{, and eighth} polytope of the predicted tube indicate vertices at which the \change{state} constraints are active. \change{RIGHT: The optimal vertex control inputs (red dots) of the optimal invariant polytope (blue shaded set). Notice that the optimal invariant set is not entirely simple---it has only $10 < 16$ isolated vertices.}}

\end{figure*}

Our first example considers a system with 
\[
A = \frac{1}{5} \left(
\begin{array}{rr}
5 & 1 \\
-1 & 4
\end{array}
\right), \ \ B = \frac{1}{5} \left(
\begin{array}{c}
0 \\
1
\end{array}
\right),
\ \ \text{and} \ \ C = \frac{1}{5} \left(
\begin{array}{rr}
1 & 0 \\
0 & 1
\end{array}
\right),
\]
where $[A, B] = [\overline A_1, \overline B_1]$ is constant, and
\begin{eqnarray}
\mathbb X &=& [-10,6.8] \times [-4.8,10], \ \ \mathbb U = [-10,10], \notag \\[0.1cm]
\text{and} \quad \mathbb W &=& [-0.5,0.5] \times [-2,2] \; . \notag
\end{eqnarray}
We start with a regular polytope with $m \geq 3$ facets and choose $Y$ and $\sigma$ as in Remark~\ref{rem::2D} with 
\begin{equation*}
\varphi_1 = 0, \ \varphi_2 =\frac{2\pi}{m}, \ \varphi_{3} = \frac{4\pi}{m} ,\ldots, \ \varphi_{m}=\frac{2(m-1)\pi}{m}\;. 
\end{equation*}
Here, $P(Y,\sigma)$ is regular and entirely simple, but it is neither contractive nor invariant. However, by Remark~\ref{rem::Y}, specifically, by solving~\eqref{eq::feasibilityProblem} for $\beta = 0.95$, $\beta$-contractive
polytopes can be found for any $m \geq 6$. Notice that this leads to polytopes with $\overline m = m$ vertices. An optimal vertex configuration domain can be found by computing $E$ as explained in Remark~\ref{rem::2D}.

\bigskip
\noindent
Let us choose $\ell$ as in~\eqref{ex::eq::ell} with $\mathsf{Q} = \mathsf{S} = \mathbb{1}$,  $\mathsf{R} = 1$, and $\mathsf{T} = 0.01$. Additionally, the prediction horizon is set to $N = 50$.
\color{black}
\change{The left part of} Figure~\ref{fig::prediction} shows the predicted tube obtained by solving~\eqref{eq::qp} for the state measurement
\[
\change{x}_0 = \hat x = [ \ \change{4, \ 8} \ ]^\tr
\]
and $m = \overline{m} = 16$. The shape and size of the tube cross sections 
change in time as all tube parameters are optimized. This illustrates the advantage of CCTMPC compared to RTMPC. The tube converges to the optimal invariant polytope $P(Y,y_\mathrm{s})$, visualized as a dark blue shaded set in the left part of Figure~\ref{fig::prediction}. Moreover, the corresponding optimal vertex control inputs $u_\mathrm{s}$ are visualized in the form of the red dots in the right part of Figure~\ref{fig::prediction}. The black dotted lines at these vertex control inputs show how they correspond to the vertices of the optimal invariant polytope $P(Y,y_\mathrm{s})$---colored in blue. Notice that it is impossible to interpolate all vertex control inputs with one hyperplane. The gray shaded areas correspond, however, to one out of many possible continuous piecewise affine control laws that interpolate all vertex inputs as well as the central input, $\overline u = 0$ (visualized as green dot).

\begin{figure}[h!]
\begin{center}
\includegraphics[scale=0.2]{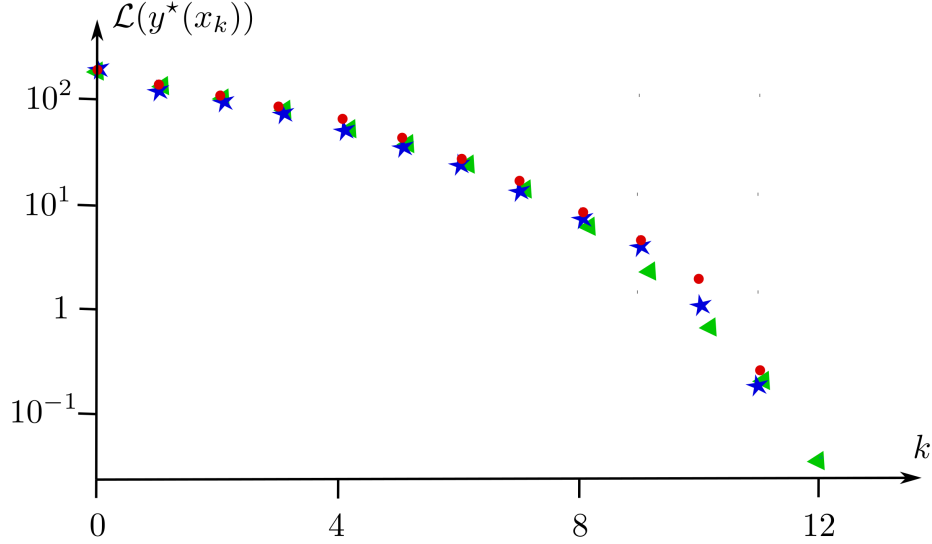}
\end{center}
\caption{Lyapunov function $\mathcal L$ along the closed-loop trajectory $\change{x}_k$ for three random uncertainty scenarios $w$. In the first scenario (red dots) and the third scenario (\change{blue stars}), the closed-loop state reaches the optimal invariant set after \change{$12$} iterations, i.e. $\mathcal L(y^\star(\change{x}_k)) = 0$ for all \change{$k \geq 12$} (not visualized, as we use a logarithmic scale). For the second scenario, the invariant set is reached after \change{$13$} iterations. In all cases, the Lyapunov function is strictly descending in each iteration. As predicted by Theorem~\ref{thm::stability}, this descent property is independent of the uncertainty sequence. 
\label{fig::lyapunov}}
\end{figure}

\bigskip
\noindent
The statement of Theorem~\ref{thm::stability} can be verified numerically by running the above MPC control setup with \change{$m=16$} in a closed loop for random uncertainty scenarios. Figure~\ref{fig::lyapunov} shows an evaluation of the Lyapunov function
$
\mathcal L(y^\star(\change{x}_k))
$
along the MPC closed-loop trajectories $\change{x}_k$ for three such randomly chosen scenarios. Notice that, the specific value of $\mathcal L$---and, thus, the actual control performance---depends on the uncertainty. However, independently of how this uncertainty is chosen, $\mathcal L$ is always strictly monotonically descending until it reaches $0$, as confirmed in all our tests.

\begin{figure*}[t]
\begin{center}
\includegraphics[width=0.99\textwidth]{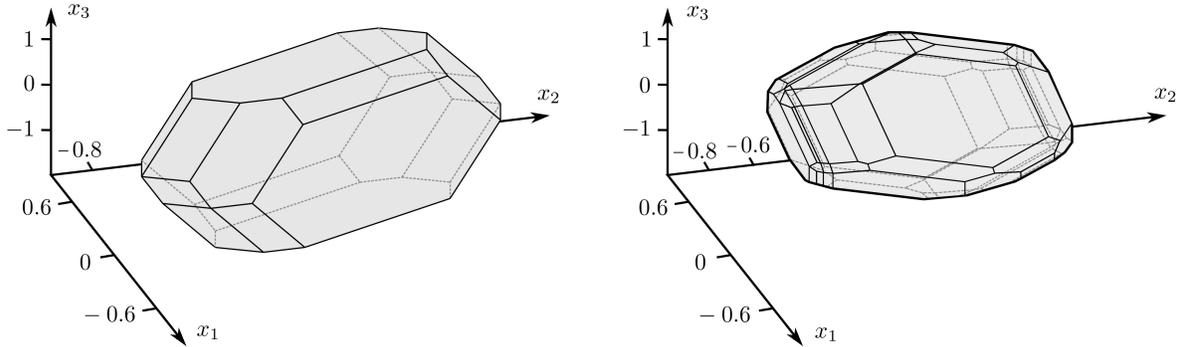}
\end{center}
\caption{\label{fig::3Dcomparison} \change{Optimal invariant polytopes for the system from Example~2. The left polytope has been found for a matrix $Y$ with $m=26$ rows and an associated vertex configuration domain that allows us to represent polytopes with up to $\overline m = 48$ vertices, while the right polytope is obtained for a setting with $m=124$ facet directions and $\overline m = 336$ vertices. Its volume is approximately equal to one third of the volume of the left polytope.}}
 \end{figure*}

\subsection{\change{Example 2}}

This section discusses an example in $\mathbb{R}^{3}$, given by
\[
A = \frac{1}{2} \left(
\begin{array}{rrr}
2 & 1 & 0 \\
0 & 2 & 1 \\
0 & 0 & 2
\end{array}
\right), \ \ B = \frac{1}{2} \left(
\begin{array}{c}
1 \\
1 \\
1
\end{array}
\right),
\ \ C = \frac{1}{2} \left(
\begin{array}{rrr}
1 & 0 & 0 \\
0 & 1 & 0 \\
0 & 0 & 1
\end{array}
\right),
\]
again with $[A, B] = [\overline A_1, \overline B_1]$ being given, as well as
\begin{eqnarray}
\mathbb X &=& [-5,5]^3, \ \ \mathbb U = [-10,10], \quad
\text{and} \quad \mathbb W = \left[ -\frac{1}{5},\frac{1}{5} \right]^3 \; . \notag
\end{eqnarray}
We choose $\ell$ as in~\eqref{ex::eq::ell} with $\mathsf{Q} = \mathbb{1}$, $ \mathsf{S} = 7 \cdot \mathbb{1}$, $\mathsf{R} = 1$, and $\mathsf{T} = 0.01$. The matrix $Y$ has row vectors of the form
\begin{align}
\label{eq::ijk}
\frac{\left[ i, \ j, \ k \right]}{\sqrt{i^2+j^2+k^2}} \; ,
\end{align}
where $(i,j,k) \in \{ -1,0,1 \}^3 \setminus \{ (0,0,0) \}$. We initially set $\sigma = [1,1,\ldots,1]^\tr$ in order to compute a vertex configuration domain and a contractive polytope as explained in Remark~\ref{rem::Y}. This leads to a contractive polytope with
\[
m = 3^3-1 = 26
\]
facets and $\overline m = 48$ vertices. Its reduced configuration matrix $E$ has $48$ rows and $168$ non-zero coefficients. The left part of Figure~\ref{fig::3Dcomparison} shows a $3$-dimensional visualization of the corresponding optimal invariant set. To illustrate how this result can be improved by increasing the number of facet directions and vertices, the right part of Figure~\ref{fig::3Dcomparison} shows a less conservative optimal invariant set constructed with $m = 5^3 - 1 = 124$ facet directions and $\overline m = 336$ vertices. Its reduced configuration matrix $E$ has $290$ rows and $1031$ non-zero coefficients. Here, the rows of $Y$ have also been constructed as in~\eqref{eq::ijk} but with $(i,j,k) \in \{ -2,-1,0,1,2 \}^3 \setminus \{ (0,0,0,0,0) \}$.

Numerical experiments indicate that ETMPC can be used to generate similar polytopes for both choices of $Y$. These are, however, at least $1 \%$ sub-optimal for the best possible affine control law (in terms of stage cost value) that we managed to compute by exhaustive search. Without optimizing the affine feedback law, ETMPC is usually more than $1 \%$ sub-optimal. 
The same trend can be observed in closed-loop simulations, which is, for the sake of brevity, not further discussed at this point.

\subsection{Example 3}
\label{sec::ex3}
Our last example shows that CCTMPC can outperform both DAFMPC and FPTMPC in terms of conservatism and run-time complexity. We consider the control system 
\begin{align}
\label{eq::system3}
\begin{array}{rcl}
x_1^+ &=& \frac{1}{2} x_3 - \frac{1}{2} x_4 \\[0.16cm]
x_2^+ &=& \frac{1}{2} x_3 + \frac{1}{2} x_4 + u \\[0.16cm]
x_3^+ &=& x_4 \\[0.16cm]
x_4^+ &=& w
\end{array} \quad \text{with} \quad \left\{ \begin{array}{rcl}
\mathbb X &=& \mathbb R^4 \\[0.16cm]
\mathbb U &=& [0,1] \\[0.16cm]
\mathbb W &=& [-1,1] \; .
\end{array} \right.
\end{align}
Here, $x \in \mathbb R^4$ denotes the  state at a given time and $x^+$ the corresponding 
successor state. Sub-indices denote the components of the states (not a time index). Our goal is to minimize the least-squares distance of all states to $0$, for example, by setting $\ell(y,u) = y^2$, but other choices of $\ell$ are possible, too. It is not difficult to see that an optimal feedback law is in this case given by
\begin{align}
\mu(x) \ = \ \left\{ 
\begin{array}{ll}
0 & \text{if} \ \ \  x_3 + x_4 > 0 \\[0.16cm]
-\frac{x_3+x_4}{2} \ \ & \text{if} \ \ \ x_3 + x_4 \in [-2,0] \\[0.16cm]
1 & \text{otherwise} \; ,
\end{array} \right.
\label{eq::muopt}
\end{align}
because only $x_2$ depends on $u$. Next, we introduce the template matrix
\[
Y = \left(
\begin{array}{rrrrrr}
0 & 0 & 1 & 0 & -2 & 0 \\[0.16cm]
0 & 0 & -1 & 2 & 2 & -1 \\[0.16cm]
0 & 0 & 1 & -1 & -1 & 0 \\[0.16cm]
1 & -1 & 0 & 0 & 0 & 0
\end{array}
\right)^\tr \in \mathbb R^{6 \times 4} \; .
\]
We start with $\sigma = (1,1,1,1,1,1)^\tr$ and use the procedure from Remark~\ref{rem::Y} to compute a contractive polytope and its associated reduced vertex configuration template,\footnote{\change{For this particular example, the vertex configuration constraint is not restrictive: there exists for every $y' \in \mathbb R^6$ a parameter $y \in \mathbb R^6$ with $E y \leq 0$ such that $P(Y,y')=P(Y,y)$. The proof of this statement is left as an exercise to the reader.}}
\[
E = \left(
\begin{array}{rrrrrr}
-1 &  -1 & 0 & 0 & 0 & 0 \\[0.16cm]
0 & 0 & -2 & -1 & -1  & -2
\end{array}
\right) \; .
\]
If we minimize $\ell(y,u) = y^2$, the optimal invariant polytope $P(Y,y_\mathrm{s})$ is given by
\[
y_\mathrm{s} = \left( \ 1, \ 1, \ 0, \ 1, \ 1, \ 0 \ \right)^\tr \; .
\]
Notice that $P(Y,y_\mathrm{s})$ is a four dimensional polytope with $6$ facets and $8$ vertices. Its vertex control inputs $u_\mathrm{s} \in \mathbb R^8$ are unique and given by
\[
\forall i \in \{ 1,2,\ldots,8\}, \qquad (u_\mathrm{s})_i = \mu(V_i y_\mathrm{s}),
\]
where $\mu$ denotes the optimal control law in~\eqref{eq::muopt}. In other words, the proposed CCTMPC method is able to find the optimal nonlinear control law. As such, for this example and this choice of $Y$ one may state that CCTMPC finds the best possible invariant polytope.

\bigskip
\noindent
Our next goal is to analyze an FPTMPC controller for~\eqref{eq::system3}. Notice that the matrix $A$ of~\eqref{eq::system3} is nil-potent, $A^3 = 0$. Thus, we may set the prediction horizon of the controller to $N=3$ (without loss of generality), as the system does anyhow neither remember uncertainties nor control inputs for longer than $3$ time steps. Next, let us work out the fully parameterized extreme partial state trajectories. They are given by 
\begin{eqnarray}
\tilde x^{0,3} &\ = \ & \left( \ 0, \ \ \tilde u^{0,2}, \ \ 0, \ \ 0 \
\right)^\tr, \notag \\[0.16cm]
\tilde x^{1,1,3} & \ = \ & \left( \ 0.5, \ \ 0.5 + \tilde u^{1,1,2}, \ \ 0, \ \ 0 \
\right)^\tr, \notag \\[0.16cm]
\tilde x^{2,1,3} & \ = \ & \left( \ -0.5, \ \ -0.5 + \tilde u^{2,1,2}, \ \ 0, \ \ 0 \
\right)^\tr, \notag \\[0.16cm]
\tilde x^{1,2,3} & \ = \ & \left( \ -0.5, \ \ 0.5 + \tilde u^{1,2,2}, \ \ 1, \ \ 0 \
\right)^\tr, \notag \\[0.16cm]
\tilde x^{2,2,3} & \ = \ & \left( \ 0.5, \ \ -0.5 + \tilde u^{2,2,2}, \ \ -1, \ \ 0 \
\right)^\tr, \notag \\[0.16cm]
\tilde x^{1,3,3} & \ = \ & \left( \ 0, \ \ 0, \ \ 0, \ \ 1 \
\right)^\tr, \notag \\[0.16cm]
\text{and} \quad \tilde x^{2,3,3} & \ = \ & \left( \ 0, \ \ 0, \ \ 0, \ \ -1 \ 
\right)^\tr \; ,
\end{eqnarray}
where we use the same index convention as in~\cite{Rakovic2012a}. That is,  
$\tilde x^{i,k,3}$ is the extreme partial state trajectory at time $3$, obtained by exciting the system at time $k-1$ with the extreme input $\tilde w_i = (-1)^{i-1}$. The value of the measurement $\tilde x^{0,0} = \hat x$ at time $0$ is irrelevant, as the state at time $3$ is independent of $\hat x$. The corresponding inputs,
\[
\tilde u = \left( \ \tilde u^{0,2}, \ \tilde u^{1,1,2}, \ \tilde u^{2,1,2}, \ \tilde u^{1,2,2}, \ \tilde u^{2,2,2} \right)^\tr \in \mathbb R^{5},
\]
need to satisfy the control constraints for all possible extreme uncertainty scenarios; that is,
\begin{align}
\label{eq::PTU1}
& \tilde u^{0,2} + \tilde u^{1,1,2} + \tilde u^{1,2,2} \ \in \ [0,1] \\[0.16cm]
\label{eq::PTU2}
& \tilde u^{0,2} + \tilde u^{1,1,2} + \tilde u^{2,2,2} \ \in \ [0,1] \\[0.16cm]
\label{eq::PTU3}
& \tilde u^{0,2} + \tilde u^{2,1,2} + \tilde u^{1,2,2} \ \in \ [0,1] \\[0.16cm]
\label{eq::PTU4}
\text{and} \qquad & \tilde u^{0,2} + \tilde u^{2,1,2} + \tilde u^{2,2,2} \ \in \ [0,1] \; .
\end{align}
Recall that $P(Y,y_\mathrm{s})$ denotes the optimal invariant set, which has $6$ facets and $8$ vertices, as discussed above. Let us assume that FPTMPC was not conservative. In this case, it would be possible to find a parameter $\tilde u$ satisfying~\eqref{eq::PTU1}--\eqref{eq::PTU4} and ensuring that all possible extreme state trajectories are contained in 
$P(Y,y_{\rm s})$; that is,
\begin{align}
\label{eq::PTS1}
& \tilde x^{0,3} + \tilde x^{1,1,3} + \tilde x^{1,2,3}  + \tilde x^{i,3,3} \in P(Y,y_\mathrm{s}), \\[0.16cm]
\label{eq::PTS2}
& \tilde x^{0,3} + \tilde x^{1,1,3} + \tilde x^{2,2,3}  + \tilde x^{i,3,3} \in P(Y,y_\mathrm{s}), \\[0.16cm]
\label{eq::PTS3}
& \tilde x^{0,3} + \tilde x^{2,1,3} + \tilde x^{1,2,3}  + \tilde x^{i,3,3} \in P(Y,y_\mathrm{s}), \\[0.16cm]
\label{eq::PTS4}
& \tilde x^{0,3} + \tilde x^{2,1,3} + \tilde x^{2,2,3}  + \tilde x^{i,3,3} \in P(Y,y_\mathrm{s})
\end{align}
for all $i \in \{1,2\}$ ($2^3$ conditions in total). A closer inspection reveals that Conditions~\eqref{eq::PTS1} and~\eqref{eq::PTU1} necessarily imply Condition~\eqref{eq::PTSU1}, as stated below. Analogous necessary conditions of feasibility are given by
\begin{align}
\label{eq::PTSU1}
\tilde u^{0,2} + \tilde u^{1,1,2} + \tilde u^{1,2,2} & \ \overset{\eqref{eq::PTU1},\eqref{eq::PTS1}}{=} \ 0, \\[0.16cm]
\label{eq::PTSU2}
\tilde u^{0,2} + \tilde u^{1,1,2} + \tilde u^{2,2,2} & \ \overset{\eqref{eq::PTU2},\eqref{eq::PTS2}}{=} \ 0, \\[0.16cm]
\label{eq::PTSU3}
\tilde u^{0,2} + \tilde u^{2,1,2} + \tilde u^{1,2,2} & \ \overset{\eqref{eq::PTU3},\eqref{eq::PTS3}}{=} \ 0, \\[0.16cm]
\text{and} \quad \label{eq::PTSU4}
\tilde u^{0,2} + \tilde u^{2,1,2} + \tilde u^{2,2,2} & \ \overset{\eqref{eq::PTU4},\eqref{eq::PTS4}}{=} \ 1 \; .
\end{align}
Subtracting~\eqref{eq::PTSU2} from~\eqref{eq::PTSU1} yields $\tilde u^{1,2,2} = \tilde u^{2,2,2}$. By substituting this equation in~\eqref{eq::PTSU3} and subtracting it from~\eqref{eq::PTSU4} we find that $1 = 0$. Clearly, this is a contradiction, which implies that FPTMPC is conservative and this result is independent of how one chooses the prediction horizon $N \geq 3$ and the initial measurement $\hat x$.

\bigskip
\noindent
In order to additionally explain why CCTMPC is for the current example less conservative than DAFMPC, it is helpful to note that the optimal feedback law~\eqref{eq::muopt} is equivalent to a nonlinear disturbance feedback control law of the form
\begin{align}
\widetilde \mu(w^{-},w^{--}) \ = \ \left\{ 
\begin{array}{ll}
0 & \text{if} \ \  w^{-}+w^{--} > 0 \\[0.16cm]
-\frac{w^{-}+w^{--}}{2} \ & \text{if} \ \ w^{-}+w^{--} \in [-2,0] \\[0.16cm]
1 & \text{otherwise} \; ,
\end{array} \right. \notag
\end{align}
where $w^-$ denotes the last and $w^{--}$ the disturbance from two time steps ago. By evaluating this feedback law at the extreme scenarios,
\begin{align}
\widetilde \mu(-1,1) = 0, \ \ & \ \widetilde \mu(1,1) = 0, \notag \\[0.16cm]
\widetilde \mu(-1,-1) = 1, \ \ & \ \widetilde \mu(1,-1) = 0, \notag
\end{align}
one finds that it is impossible to interpolate these $4$ uniquely optimal function values with a single affine function. This observation implies that, for this example, CCTMPC is strictly less conservative than DAFMPC.

\section{Conclusions}
\label{sec::conclusions}
This paper has presented a novel class of Tube MPC controllers \change{for linear systems with additive and multiplicative uncertainty}. The corresponding technical developments built upon Theorem~\ref{thm::convhull}, which \change{features a variant of the Gauss-Bonnet theorem in order to simultaneously parameterize the facets and vertices of configuration-constrained polytopes. The relevance of this geometrical construction is that it enables us to freely optimize configuration-constrained robust forward invariant tubes and their associated vertex control laws via the convex optimization problem~\eqref{eq::qp}.}
Conditions under which \change{the resulting CCTMPC controller is asymptotically stable and recursively feasible} are established in Theorems~\ref{thm::Lyapunov} and~\ref{thm::stability}.

\bigskip
\noindent
CCTMPC is never more conservative than RTMPC, HTMPC, and---in many cases, for instance, for all systems with $2$ states---also ETMPC. Moreover, we have constructed examples for which CCTMPC is less conservative than FTPMPC and DAFMPC. Additionally, for sufficiently long time horizons, CCTMPC can be guaranteed to be computationally less demanding than FTPMPC and DAFMPC. And, finally, a unique advantage of CCTMPC compared to existing robust MPC schemes is that CCTMPC can naturally take additive and multiplicative uncertainties into account.

\newpage

\bibliographystyle{plain}
\bibliography{references}

\appendix

\section{Proof of Theorem~\ref{thm::regularTemplates} \change{and Corollary~\ref{cor::Prevalent}}}
\label{sec::proofA}
We first show that the second statement of Theorem~\ref{thm::regularTemplates} implies the first statement. The reverse implication is established in Part~II, \change{which will also imply the statement of Corollary~\ref{cor::Prevalent}.}

\textit{Part I.} Let $\sigma \in \mathbb Y_{\mathcal I}$ be such that $P(Y,\sigma)$ is an entirely simple polyhedron and let $I \in \mathcal I$ be a given index set. First, we 
show
\[
J = \left\{ \ j \notin I \ \middle| \ 0 = \min_{x \in \mathcal F_I(\sigma)} \, Y_j x - \sigma_j \ \right\} = \varnothing\;.
\]
Notice that $\mathcal F_{I \cup J}(\sigma) = \mathcal F_I(\sigma)$ holds by construction and both sets are non-empty, since $\sigma \in \mathbb Y_{\mathcal I}$ and $I \in \mathcal I$. Moreover, there exists a point $x^\star \in \mathcal F_{I \cup J}(\sigma)$ that satisfies
\begin{align}
\label{eq::intJI}
\forall k \notin J \cup I, \quad Y_k x^\star - y_k \; < \; 0 \; ,
\end{align}
which directly follows from the definition of $J$. Next, let
\begin{align}
\mathcal T_{I \cup J}(\sigma) &= \left\{ \ x \in \mathbb R^n \ \middle| \ Y_{I \cup J} x = \sigma_{I \cup J} \ \right\} \\[0.16cm]
\text{and} \qquad \quad \mathcal T_{I}(\sigma) &= \left\{ \ x \in \mathbb R^n \ \middle| \ Y_{I} x = \sigma_{I} \ \right\}
\end{align}
denote the linear subspaces of $\mathbb R^n$ corresponding to $I \cup J$ and $I$. Since $P(Y,\sigma)$ is entirely simple,~\eqref{eq::intJI} implies that $\mathcal F_{I \cup J}(\sigma)$ has non-empty interior in $\mathcal T_{I \cup J}(\sigma)$, implying that the dimension of the subspace $\mathcal T_{I \cup J}(\sigma)$ and the face $\mathcal F_{I \cup J}(\sigma)$ must coincide,
\[
\mathrm{dim}(\mathcal F_{I \cup J}(\sigma)) = \mathrm{dim}(\mathcal T_{I \cup J}(\sigma)) \; .
\]
Consequently, since $\mathcal F_{I \cup J}(\sigma) = \mathcal F_{I}(\sigma)$, we have
\begin{align}
n - \rank( Y_{I} ) &= \mathrm{dim}( \mathcal T_{I}(\sigma) ) \, \geq \, \mathrm{dim}( \mathcal F_{I}(\sigma) ) \, = \, \mathrm{dim}( \mathcal F_{I \cup J}(\sigma) ) \notag \\[0.16cm]
&= \mathrm{dim}(\mathcal T_{I \cup J}(\sigma)) \,
\label{eq::rankYY}
=  n - \rank( Y_{I \cup J} ) \; .
\end{align}
Since $P(Y,\sigma)$ is entirely simple, the faces $\mathcal F_{I}$ and $\mathcal F_{I \cup J}$ satisfy
\begin{align}
\label{eq::rankYI}
\rank(Y_I) &= |I| \\[0.16cm]
\label{eq::rankYIJ}
\text{and} \qquad \rank(Y_{I \cup J}) &= |I \cup J| = |I| + |J| \; ,
\end{align}
where~\eqref{eq::rankYIJ} holds, as $I$ and $J$ are, by construction, disjoint. Now, we can substitute~\eqref{eq::rankYI} and~\eqref{eq::rankYIJ} in~\eqref{eq::rankYY}, which yields
\[
n - |I| \; \geq \; n - |I| - |J| \quad \Longrightarrow \quad |J| = 0 \; .
\]
Consequently, $J = \varnothing$. In other words, there exists a point $x^\star \in \mathcal F_I(\sigma)$ such that
\begin{align}
\label{eq::strictFeasFI}
\forall j \notin I, \qquad Y_j x^\star < \sigma_j \quad \text{and}  \quad Y_I x^\star = \sigma_I\; .
\end{align}
Next, let us perturb $\sigma$ by a small vector $\delta$. Since $Y_I$ has full rank, 
$x_{\delta}^\star \defeq x^\star + Y_I^\dagger \delta$ is well-defined for any perturbation
$\delta \in \mathbb R^m$. Here, $Y_I^\dagger$ denotes the right pseudo-inverse of $Y_I$. Since~\eqref{eq::strictFeasFI} holds, there exists a small $\epsilon_I > 0$ such that
\begin{align}
\label{eq::strictFeasFI2}
\forall j \notin I, \quad Y_j x_\delta^\star < \sigma_j + \delta_j \quad \text{and}  \quad Y_I x_\delta^\star = \sigma_I + \delta_I\; .
\end{align}
for all $\delta \in \mathbb R^m$ with $\| \delta \| < \epsilon_I$. Consequently, the face $\mathcal F_I(\sigma + \delta)$ is non-empty for any such small perturbation. Thus, if we set
\[
0 < \epsilon \; \defeq \; \min_{I \in \mathcal I} \; \epsilon_I \; ,
\]
it follows that $\sigma + \delta \in \mathbb Y_{\mathcal I}$ for all $\delta \in \mathbb R^m$ with $\| \delta \| \leq \epsilon$. But this means that $\mathbb Y_{\mathcal I}$ has a non-empty interior in $\mathbb R^m$ and, consequently, it is a regular configuration domain. This concludes the first part of the proof.

\textit{Part II.} Let us introduce the notation
\[
\mathcal R_I \; \defeq \; \left\{ \ y \in \mathbb R^m \ \middle| \
\begin{array}{l}
\mathcal F_I(y) \neq \varnothing,\ \rank(Y_I) < |I|
\end{array}
\ \right\}
\]
implying that if $R_I \neq \varnothing$, the matrix $Y_I$ is degenerate in the sense that it has at least one degenerate row, which is not linearly independent of its other rows. In other words, there exists a linear subspace $T_I \subseteq \mathbb R^m$ of $\mathbb R^m$ (namely the one belonging to the degenerate row of $Y_I$) such that
\[
\mathcal R_I \subseteq T_I \quad \text{and} \quad \mathrm{dim}(T_I) \leq m-1 \; .
\]
Now, observe that the set
\[
\mathcal R \; \defeq \; \mathbb Y \setminus \left(  \bigcup_{I \subseteq \{ 1, \ldots, m \} } \mathcal R_I \right)
\]
corresponds---by definition---to the set of parameters $y$ for which $P(Y,y)$ is an entirely simple polyhedron. Because $\mathbb Y$ has non-empty interior in $\mathbb R^{m}$ (see Corollary~\ref{cor::nonempty}) and because the sets $\mathcal R_I$ are contained in the not full-dimensional subspaces $T_I$ of $\mathbb R^m$, the set $\mathcal R$ is dense in $\mathbb Y$. Thus, if $\mathbb Y_{\mathcal I}$ is a regular configuration domain, the set
$\mathcal R \cap \mathbb Y_{\mathcal I}$ is non-empty. That is, we can find a point $\sigma \in \mathbb Y_{\mathcal I}$ such that $P(Y,\sigma)$ is an entirely simple polyhedron. \change{Notice that the same argument implies that the statement of Corollary~\ref{cor::Prevalent} holds, as the complement of $\mathcal R$ in $\mathbb Y$ has Lebesgue measure zero.}

\section{Proof of Corollary~\ref{cor::LocalStability}}
\label{app::B}
By definition, $P(Y,y)$ is locally configuration stable if and only if $y \in \mathrm{int}( \mathbb Y_{\mathcal C(y)} )$. Therefore, the "if" part of the corollary follows as in Part I of the above proof of Theorem~\ref{thm::regularTemplates} by replacing $\sigma$ by $y$. Concerning the "only if" part, note that if $\mathcal F_I(y)$ is non-empty, the equation $Y_I x_I = y_I$ is feasible. If $Y_I$ has a degenerate row, there exist arbitrarily small perturbations of $y$ for which this equation becomes infeasible and the face configuration changes. Thus, by using the definition of $\mathcal R$ from Part II of the above proof of Theorem~\ref{thm::regularTemplates}, we find that
$\mathcal R \cap \mathbb Y_{\mathcal C(y)} \supseteq \mathrm{int}( \mathbb Y_{\mathcal C(y)} )
$. Since $y \in \mathrm{int}( \mathbb Y_{\mathcal C(y)} )$, the definition of $\mathcal{R}$ implies that $P(Y,y)$ is entirely simple.
\end{document}